\documentclass[11pt,twoside]{article}
\usepackage{graphicx}
\usepackage{fancyhdr}
\usepackage[cmex10]{amsmath}
\usepackage{amsmath, amssymb, epsfig, amsthm}
\usepackage{bm}
\usepackage{mathrsfs}
\usepackage{color}
\usepackage{algorithm}
\usepackage{algorithmic}
\usepackage{array}
\usepackage[font=footnotesize]{subfig}

\newcommand{\Ik}{{I^{(k)}}} 
\newcommand{\Jk}{{J^{(k)}}} 
\newcommand{\IkcapT}{{I^{(k)}\cap T}} 
\newcommand{\JkcapT}{{J^{(k)}\cap T}} 
\newcommand{\IkminusT}{{I^{(k)}-T}} 
\newcommand{\JkminusT}{{J^{(k)}-T}} 
\newcommand{\cT}{{c_T}} 
\newcommand{\cIk}{{c_\Ik}} 
\newcommand{\cJk}{{c_\Jk}} 
\newcommand{\chatT}{{\hat c_T}}
\newcommand{\chatIk}{{\hat{c}_\Ik}}  
\newcommand{\chatJk}{{\hat{c}_\Jk}} 
\newcommand{\chatk}{{\hat{c}_{(k)}}} 
\newcommand{\cIkcapT}{{c_\IkcapT}} 
\newcommand{\cJkcapT}{{c_\JkcapT}}

\newcommand{\chatIkcapT}{{\hat c_\IkcapT}} 
\newcommand{\chatJkcapT}{{\hat c_\JkcapT}} 
\newcommand{\chatIkminusT}{{\hat c_\IkminusT}} 
\newcommand{\chatJkminusT}{{\hat c_\JkminusT}} 
\newcommand{\PhiT}{{\Phi_T}} 
\newcommand{\PhistarT}{{\Phi^*_T}} 
\newcommand{\PhiIk}{{\Phi_\Ik}} 
\newcommand{\PhiJk}{{\Phi_\Jk}} 
 
\newcommand{\PhistarJk}{{\Phi^*_\Jk}} 
\newcommand{\PhiIkpinv}{{\Phi^\dag_\Ik}} 

\newcommand{\PhiIkcapT}{{\Phi_\IkcapT}} 
\newcommand{\PhiJkcapT}{{\Phi_\JkcapT}} 
 
\newcommand{\PhistarJkcapT}{{\Phi^*_\JkcapT}}

\newcommand{\PhistarJkminusT}{{\Phi^*_\JkminusT}}
\DeclareMathOperator{\supp}{supp} 
\theoremstyle{plain}
\newtheorem{theorem}{Theorem}[section]

\theoremstyle{definition}

\theoremstyle{remark}

\numberwithin{equation}{section}

\usepackage{fullpage}

\begin{document}
\title{Guaranteed sparse signal recovery with highly coherent sensing matrices}

\author{Guangliang Chen \\ \small Department of Mathematics, Claremont McKenna College \\
\small Claremont, California 91711, USA  \\ \small gchen@cmc.edu\\
\\
 Atul Divekar \\ \small Alcatel-Lucent,\\
\small Naperville, Illinois 60563, USA \\ \small atul.divekar@alcatel-lucent.com \\
\\
 Deanna Needell \\ \small Department of Mathematics, Claremont McKenna College \\
\small Claremont, California 91711, USA  \\ \small dneedell@cmc.edu }

\newcommand{\vnorm}[1]{||#1||}
\newcommand{\vns}[1]{{\vnorm{#1}^{2}}}
\newcommand{\IOmega}[1]{{I_{1} \cap \Omega_{#1}}}
\newcommand{\IRminusOmega}[1]{{I_{1}\cap R \backslash \Omega_{#1}}}
\newcommand{\IOmegaminusR}[1]{{I_{1}\cap \Omega_{#1}\backslash R}}
\newcommand{\GOmega}[1]{{G \cap \Omega_{#1}}}
\newcommand{\RminusOmega}[1]{{R \backslash \Omega_{#1}}}
\newcommand{\GRminusOmega}[1]{{G \cap \RminusOmega{#1}}}
\newcommand{\OmegaminusR}[1]{{\Omega_{#1}\backslash R}}
\newcommand{\GOmegaminusR}[1]{{G \cap \OmegaminusR{#1}}}

\maketitle
\begin{abstract}
Compressive sensing is a methodology for the reconstruction of sparse 
or compressible signals using far fewer samples than required by the
Nyquist criterion.  However, many of the results in compressive sensing
concern random sampling matrices such as Gaussian and Bernoulli matrices.
In common physically feasible signal acquisition and reconstruction 
scenarios such as super-resolution of images, the sensing matrix has
a non-random structure with highly correlated columns. Here we present a compressive sensing recovery algorithm, called \emph{Partial Inversion (PartInv)}, that shows better performance than existing greedy methods for random matrices, and is especially suitable for matrices that have subsets of highly correlated columns. We provide theoretical justification as well as empirical comparisons.
\end{abstract}

\section{Introduction}
 
Consider the problem of image super-resolution, where one or more low-resolution images of a scene are used to synthesize a single image of higher resolution. If multiple images are used, they are commonly assumed to be subpixel-shifted and downsampled versions of the original high resolution image that is to be reconstructed~\cite{superres_survey}. Alternatively, super-resolution from a single low resolution image using a dictionary of image patches and compressive sensing recovery has been proposed in \cite{wright_superres}. The relationship between the available low resolution and desired high resolution images is commonly modeled by a linear filtering and downsampling operation. Suppose that we wish to reconstruct a size $N\times N$ high resolution image from a lower resolution image, for example of size $\frac{N}{2}\times \frac{N}{2}$, or smaller.  Let $x$ and $y$ represent the vectorized high and low resolution images respectively. We model the formation of $y$ from $x$ by the equation $y=SHx + \eta$ where $\eta$ is the sensor noise, $S$ is a downsampling matrix of size $\frac{N^2}{4}$ by $N^{2}$, and $H$ is an $N^{2}$ by $N^{2}$ matrix that represents the filtering (antialiasing) operation. In order to consider super-resolution as a compressive sensing recovery problem we write $x=\Psi c$ where $\Psi$ is a sparsifying basis for the class of images under consideration and $c$ is the coefficient vector corresponding to image $x$ with respect to the basis $\Psi$.  Typically, one assumes the signal $x$ has a sparse representation, meaning that the coefficient vector $c$ has a small number of non-zero coefficients: $\|c\|_0 := |\supp(c)| = K \ll N$.  We call such vectors $K$-sparse, and for approximately sparse vectors we write $x_K$ to denote the best $K$-sparse approximation of $x$.
In the simplest case, $\Psi$ is an $N^{2}\times N^{2}$ orthogonal matrix, but can also be generalized to an overcomplete dictionary.
Here we have 
$
y = SH\Psi c + \eta = \Phi c + \eta, 
$
where $\Phi=SH\Psi$ is the sampling matrix. 

 Most of the work in the compressive sensing literature assumes $\Phi$ to be a random matrix, such as a partial DFT or one drawn from a Gaussian or Bernoulli distribution. However, in this scenario
the matrix is not random, but instead has correlated columns whose structure may impair conventional compressive sensing recovery. Here $H$ may not be a perfect low pass filter, so that it is possible for $\Phi=SH\Psi$ to preserve enough high frequency information for recovery to be possible; $SH$ and $\Psi$ have sufficient incoherency to allow $c$ to be recovered with acceptable error.

 Compressed sensing provides techniques for stable sparse recovery~\cite{candes2006robust,candes2006stable,donoho1989uncertainty}, but results for coherent sensing matrices have been limited~\cite{candes2011compressed,fannjiang2012coherence,candes2012towards}.
 
 {\bf Contribution.}  In this work we present an algorithm called the Partial Inversion (PartInv) method.
PartInv eliminates a source of noise in the proxy used by existing greedy algorithms such as CoSaMP to identify the nonzero coefficients, which gives it better performance, since there is always some correlation among the columns of a rectangular sensing matrix (i.e. one with more columns than rows). PartInv works especially well when the sensing matrix has subsets of heavily correlated columns, which is seen in the super-resolution problem. We present theoretical justification for PartInv in the setting of coherent sensing matrices and provide rigorous results for the setting of more incoherent matrices. Our experiments confirm that PartInv yields improved recovery for both cases -when matrices are coherent (the super-resolution case) and when they are more incoherent. We believe that these results can be extended to the coherent case, and leave such an analysis for future work.
 
{\bf Organization.} The correlation structure in the columns of the sensing matrix is first described in Section~\ref{sec:corr}.  In Section~\ref{sec:partinv} we present an algorithm which tackles the correlation in the sensing matrix to recover the signal, and provide theoretical justification for it.  
Section~\ref{sec:exp} contains experimental results for the algorithm and comparisons to other existing methods, as well as a discussion of the particular case involving wavelet sparsity.  We conclude in Section~\ref{sec:conc} and include the proofs of our main results in the appendix.

\section{Correlation Structure}\label{sec:corr}

Typical examples of sparsifying bases $\Psi$ for images are wavelets and blockwise discrete cosine transform bases. Images exhibit correlation at each scale: neighboring pixels are heavily correlated except across edges, local averages of neighboring blocks are heavily correlated except across edges, and so on. This makes wavelet-like bases, which have locally restricted atoms, suitable for sparsifying the image. For the super-resolution setting with the low resolution image of size $\frac{N}{2}\times \frac{N}{2}$, the rows of $SH$ consist of shifted versions of the filtering kernel with shifts of 2 horizontally and vertically. Due to the localized nature of wavelet bases, we expect columns of $\Phi$ that correspond to spatially distant bases in $\Psi$ to have little correlation. If $\Psi$ is a tree structured orthogonal wavelet basis matrix, columns of $\Psi$ that overlap spatially are orthogonal, however when filtered by $H$, they result in significant correlation. Then we expect columns in $\Phi$ to show significant correlation in tree structured patterns.

We illustrate this with an example. For simplicity we consider only one-dimensional signals,
though the discussion is equally valid for images. Suppose that $\Psi$ is a $256\times 256$ matrix whose columns consist of the length 256 Haar basis vectors, and $SH$ is a $128\times 256$ matrix obtained by shifting the filter kernel $h=\{0.1,0.2,0.4,0.2,0.1\}$ by two from one row to the next. $SH$ represents the filtering and downsampling operation that generates the low resolution signal $y=SHx$ from the length 256 signal $x$. Then $\Phi=SH\Psi$ is the sampling matrix. 

\begin{figure}[ht]
\vspace{-1mm}
\centering
  \includegraphics[height=5.0cm,width=5.0cm]{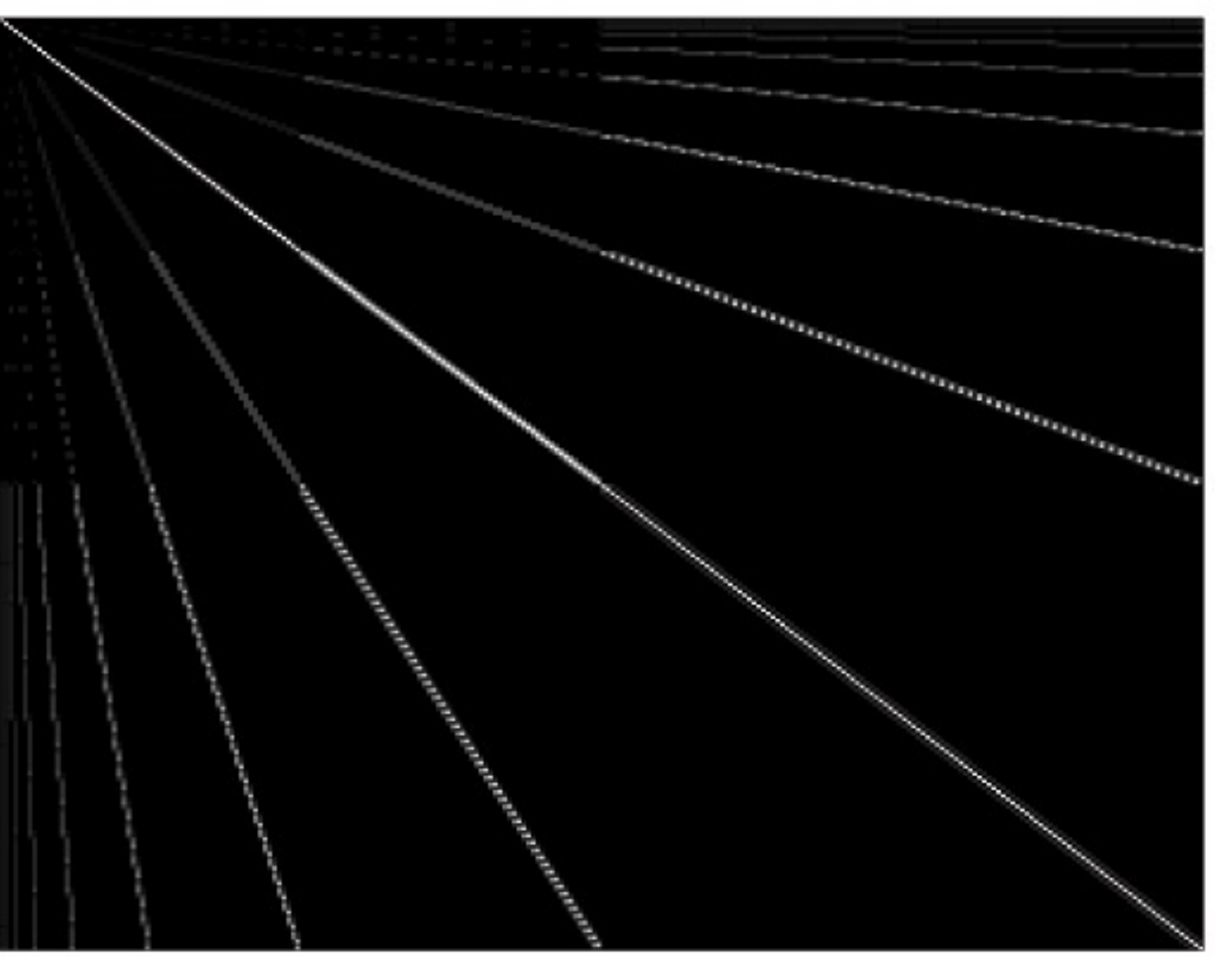} 
  \vspace{-1mm}
  \caption{Absolute values of $\Phi^{*}\Phi$.\label{corr_struct}}
  \vspace{-1mm}
\end{figure}

Fig. \ref{corr_struct} 
shows the absolute values of the correlation matrix $C=\Phi^{*}\Phi$ (here and throughout $A^*$ denotes the adjoint of $A$). This shows that only a small number of pairs of columns of $\Phi$ are strongly correlated to each other. Each filtered wavelet basis is correlated with other spatially overlapping bases at coarser and finer scale and in the immediate neighborhood, but has no correlation with spatially distant bases.

One of the central results in compressive sensing is that if matrix $\Phi$ exhibits a property called the Restricted Isometry Property (RIP)~\cite{decode_tao,RV08:sparse},
convex optimization can recover the sparse signal 
exactly~\cite{candes_tao,candes} via the following program 
\begin{equation} \label{eq:l1}
\text{min } \vnorm{c}_{1}  \qquad 
 \text{such that} \quad y= \Phi c
\end{equation}
where $ \vnorm{c}_{1}=\sum_i |c_i| $ is the $\ell_1$-norm of $c$.
However, the sampling matrix $\Phi=S H\Psi$ described above does not obey the RIP and these results are not readily applicable. 
We propose a simple modification to CoSaMP\cite{CoSaMP}, called
\emph{Partial Inversion (PartInv)} and described by Algorithm~\ref{PartInv}, to deal with correlated columns. As we shall see below, it provides an improvement in recovery performance over CoSaMP in this setting.

\section{Partial Inversion}\label{sec:partinv}

Consider the usual CS setting: Given a length $M$ sample vector $y=\Phi c + \eta$ where $\Phi$ is an $M\times N$ sampling matrix and $c$ a length $N$ vector with sparsity $K < M$, we wish to obtain the best $K$-sparse approximation $\hat{c}$ to $c$. At each step let $I$ be an index set and $\hat{c}_{I}$ represent an estimate of the components of $c$ corresponding to the column indices in $I$. The vector $\hat{c}$ by itself is an estimate for all the columns 
$\{1..N\}$. Let $L$ for $K\leq L<M$ be an adjustable parameter for the size of the 
set $I$. The value of $L$ that provides the best performance depends on the matrix, a detailed is given in a later section. Let $\Phi_{I}$ denote the matrix of columns from $\Phi$ corresponding to indices in the set $I$.  Similarly, $c_I$ denotes the vector $c$ with entries in the complement of $I$ set to zero.  At times we will write $c_I$ instead to be the vector in $\mathbb{C}^{|I|}$ consisting of the elements of $c$ indexed by $I$, in which the notation will be clear from context.  Let $\tilde{I}=\{1..N\} \backslash I $ denote the complement of $I$. For any full rank matrix $A$, define $A^{\dagger}=(A^{*}A)^{-1}A^{*}$.

\begin{algorithm}[htbp]
\caption{Given $y=\Phi c$, return best $K$-sparse approximation $\hat{c}$ }
\label{PartInv}
\begin{algorithmic}[1]
\STATE $\hat{c} \leftarrow \Phi^{*}y; I^{(0)} \leftarrow$ indices of the $L$-largest magnitudes of $\hat{c} ;  k \leftarrow 0$
\WHILE{Stopping condition not met} 
\STATE $\hat{c}_{I^{(k)}} \leftarrow \Phi_{I^{(k)}}^{\dagger} y$
\STATE $r \leftarrow y - \Phi_{I^{(k)}}\hat{c}_{I^{(k)}}$
\STATE $J^{(k)} \leftarrow \widetilde{I^{(k)}}$ 
\STATE $\hat{c}_{J^{(k)}} \leftarrow \Phi_{J^{(k)}}^{*}r$
\STATE $I^{(k+1)} \leftarrow$ indices of the $L$-largest magnitude components of $\hat{c}$.
\STATE $ k \leftarrow k+1$
\ENDWHILE
\end{algorithmic}
\end{algorithm}

For the noiseless case $\eta=0$, the stopping condition can be obtained by testing the 
magnitude of $r_{2}= y-\Phi \hat{c}$ at the start of each iteration.  If the set $I$ does not
vary from one iteration to the next, the algorithm cannot progress further and should be 
stopped immediately.  In practice the inversion of line 3 can be done efficiently by Richardson's algorithm (see e.g.~\cite[Sec.~7.2]{Bjo96:Numerical-Methods}).  

This algorithm demonstrates improvement relative to CoSaMP when the accurate recovery region
is considered on a plot of $\frac{K}{M}$ versus $\frac{M}{N}$. The motivation is the following (for simplicity we drop the iteration indicator $k$) :
From line 3, we obtain
\begin{align}
\label{eqn:PartInv}
\hat{c}_{I} &= \Phi_{I}^{\dagger} y =  \Phi_{I}^{\dagger} (\Phi_{{I}}c_{{I}}+\Phi_{\tilde{I}}c_{\tilde{I}})\\
            &= c_{I} + (\Phi_{I}^{*}\Phi_{I})^{-1}\Phi_{I}^{*}\Phi_{\tilde{I}}c_{\tilde{I}}.
\end{align}
In CoSaMP, $\Phi^*$ rather than $\Phi_{I}^{\dagger}$ is used to form a proxy and identify large coefficients of the signal.  In this case, the proxy
$\hat{c}$ restricted to the index set $I$ when $r=y$ satisfies
\begin{align*}
\hat{c_{I}} &= \Phi_{I}^{*}y \\
            &= \Phi_{I}^{*}\Phi_{I}c_{I} + \Phi_{I}^{*}\Phi_{\tilde{I}}c_{\tilde{I}} \\
            &= c_{I} + (\Phi_{I}^{*}\Phi_{I}-I)c_{I} + \Phi_{I}^{*}\Phi_{\tilde{I}}c_{\tilde{I}}
\end{align*}

If the index set $I$ contains several nonzero coefficients
(which we hope is true),  then $(\Phi_{I}^{*}\Phi_{I}-I)c_{I}$, which results from the mutual interference between the columns of $\Phi_{I}$, is significant and is a source of noise in $\hat{c_{I}}$. This term is eliminated in~\eqref{eqn:PartInv}. Partial inversion does add $(\Phi_{I}^{*}\Phi_{I})^{-1}$ to the remaining noise term, however, the singular values of this term can be kept from significantly amplifying the noise term by a conservative choice of $L$, the size of the index set $I$. For example, empirically we find that $L=K$ tends to be a safe choice, but the value of $L$ that produces optimal performance depends on the type of matrix, besides $K$ and $M$. 
The improved estimate $\hat{c}_{I}$ further produces an improved estimate $\hat{c}_{J^{(k)}}$, which leads to a better selection of nonzero coefficients in the next iteration.

 The expression~\eqref{eqn:PartInv} also indicates how the correlation structure 
may be used to improve recovery. The noise term $(\Phi_{I}^{*}\Phi_{I})^{-1}\Phi_{I}^{*}\Phi_{\tilde{I}}c_{\tilde{I}}$ depends upon the correlation between the sets $\Phi_{I}$ and $\Phi_{\tilde{I}}$ given by $\Phi_{I}^{*}\Phi_{\tilde{I}}$. This correlation is weak if $\Phi_{I}$ and $\Phi_{\tilde{I}}$ are sufficiently spread. 
However, the correlation is likely to remain large if $L$ is significant compared to $M$, as will be the case when $\frac{K}{M}$ is large.
 
In the wavelet case, if set $I$ contains several wavelet subsets and $\tilde{I}$ 
contains columns that are not from any of the subsets in $I$, then the correlation in the noise term is small.  By these arguments, we see that if the matrix $\Phi$ satisfies some mild assumptions, then Partial Inversion will converge to the sparse solution in a fixed number of iterations.  Experimentally we see that even with high correlations in the matrix $\Phi$, PartInv provides accurate recovery (see~Sec.~\ref{sec:exp}).  Under slightly stronger assumptions, we have the following mathematical justification, whose proof can be found in the appendix.  Here and throughout, we use $\|\cdot\|_2$ and $\|\cdot\|$ to denote the usual $\ell_2$ norm and spectral norm, respectively.

 \begin{theorem}\label{main}
 Let $c\in \mathbb{C}^N$ be a $K$-sparse vector with support set $T$ satisfying
\begin{align}\label{c_reqs}
|c_i| &\ge 3\delta \|c\|_2, \quad \forall i \in T,
\end{align}
for some fixed constant $0<\delta\le \frac{1}{3\sqrt{K}}$.
Assume that the dictionary $\Phi$ satisfies the following properties: 
\begin{align}
\|\Phi^*_{T_1} \Phi_{T_1} c_{T_1}\|_2 &\ge (1-\delta)^2 \|c_{T_1}\|_2 , &\quad \forall\  T_1\subseteq T \\
\|\Phi_I\| & \le A, \quad   &\forall \ |I| \leq L \\
\|\Phi_I^\dag \| & \le A, \quad   &\forall \ |I|\leq L \\
\|\Phi_I \Phi_I^\dag \Phi_{I^c \cap T }\| &\le {\delta}/{A}, \quad &\forall \ |I|\leq L.\label{las}
\end{align}
where $L$ is the parameter used in the Partial Inversion (PartInv) algorithm, and $1\le A<\sqrt{L}$ is another fixed constant.  Then PartInv reconstructs the signal, $\hat{c} = c$ in at most $K$ iterations.
 \end{theorem}

\begin{remarks}

{\bfseries 1.} First, we remark that the assumptions of this theorem restrict not only the sparsity of the signal, but also the distribution of its non-zero coefficients, contrary to typical results in compressive sensing.  We also comment that if~\eqref{c_reqs} does not hold, then the proof guarantees that all coefficients $c_i$ of $c$ that do satisfy that bound are still recovered.

{\bfseries 2.} We next relate the first three assumptions on $\Phi$ of Theorem~\ref{main} to the Restricted Isometry Property (RIP)~\cite{decode_tao}, which states that
\begin{equation}
(1-\delta)\|x\|_2^2 \leq \|\Phi x\|_2^2 \leq (1+\delta)\|x\|_2^2 \quad\text{for all $k$-sparse $x$}.
\end{equation}
It is now well-known that many classes of random matrices satisfy this property with high probabilty for $\delta < \epsilon$ when the number of measurements $M$ is on the order of $k\log(N)/\epsilon^2$~\cite{candes_tao,RV08:sparse}.
Note that the RIP is equivalent to asking that
\begin{equation}
\sqrt{1-\delta} \le \sigma_{\min} (\Phi_I) =\sigma_k (\Phi_I) \le \cdots \le \sigma_1 (\Phi_I) = \sigma_{\max} (\Phi_I)  \leq \sqrt{1+\delta}
\end{equation}
for all $|I|\le k$. Here, $\sigma_1(\cdot), \ldots, \sigma_k(\cdot)$ denote the $k$ sorted singular values (in decreasing order), and $\sigma_{\min}(\cdot), \sigma_{\max}(\cdot)$ denote the smallest and largest singular values. In contrast, our first three assumptions on $\Phi$ above are all relaxations of the RIP condition (in both directions). Indeed, we may rewrite them as follows:
\begin{align*}
\sigma_{\min}( \Phi_{T_1}) &\ge 1-\delta, &\quad \forall\  T_1\subseteq T \\
\sigma_{\max}(\Phi_I) & \le A, \quad   &\forall \ |I| \leq L \\
\sigma_r({\Phi_I}) & \ge 1/A \ (r=\mathrm{rank}(\Phi_I)), \quad   &\forall \ |I|\leq L.
\end{align*}
Note in particular that the last condition above says that the smallest \emph{nonzero} singular value of $\Phi_I$ is sufficiently large (i.e., at least $1/A$). The matrix  $\Phi_I$ may still have several zero singular values, which is in contrast to the RIP which requires the matrix $\Phi_I$ to have no zero singular values and its smallest singular value to be close to one.  This means that off the support, the matrix may have significant correlations.  The most restrictive condition is~\eqref{las}, which we believe can be improved.  We leave a detailed analysis of the requirement on $M$ in the correlated case for future work. 

{\bfseries 3.} 
We note that if the matrix $\Phi$ satisfies the RIP with parameter $L$, then the first three $\Phi$ assumptions of the theorem hold with constant $A=\sqrt{1+\delta}$ (see~\cite[Prop.~3.2]{Romp} and \cite[Prop.~3.2 and 3.3]{needell2009cosamp}). 
Since many classes of $M\times N$ random matrices satisfy the RIP when $M$ is on the order of $k\log(N/k)/\delta^2$~\cite{candes_tao,RV08:sparse}, for the assumptions of Theorem~\ref{main} to hold, one needs $M\approx k^2\log(N/k)$ measurements.  We expect that the squared dependence on $k$ is an artifact of the proof and can be improved.

{\bfseries 4.} Experimentally, we see that PartInv is also robust to noise, and doesn't necessarily require the columns of $\Phi$ to satisfy strong incoherence bounds.  We leave this theoretical analysis as future work, and provide experimental evidence in Section~\ref{sec:conc}.
\end{remarks}

\section{Experimental Results}\label{sec:exp}

We compare the recovery performance of Partial Inversion with CoSaMP and convex optimization~\eqref{eq:l1} for two classes of matrices: Gaussian random matrices, and matrices constructed to have highly correlated subsets of columns with low correlation across subsets. 

In the first case, we construct $M$ by $N$ matrices with $N(0,1)$ elements along with the coefficient vector $c$ containing $K$ nonzero entries taken from a $N(0,1)$ distribution. The nonzero locations are selected uniformly at random from $\{1...N\}$. The columns in each matrix are normalized to have unit $l_{2}$ norm. We set $N=256$ and vary $\delta=\frac{M}{N}$ from 0.1 to 0.9 in steps of 0.1. For each $\delta$ we vary $\rho=\frac{K}{M}$ from 0.1 to 0.9 in steps of 0.1. For each $(\delta,\rho)$ point we carry out 25 trials, and declare success if $\frac{1}{N}||c-\hat{c}||^{2}< 10^{-5}$. For PartInv we considered two cases for the size of subset $I$ : $L=K$ and $L=\max\{K,0.8M\}$. We see better performance in the $L=K$ case.  For $l_{1}$ minimization we use the $l_{1}$-magic package \cite{l1_magic}.  

\begin{figure}[ht]
\vspace{-1mm}
   \centering
   \begin{tabular}{@{\hspace{-2mm}}c@{\hspace{-3mm}}c}
   (a) \includegraphics[width=1.6in]{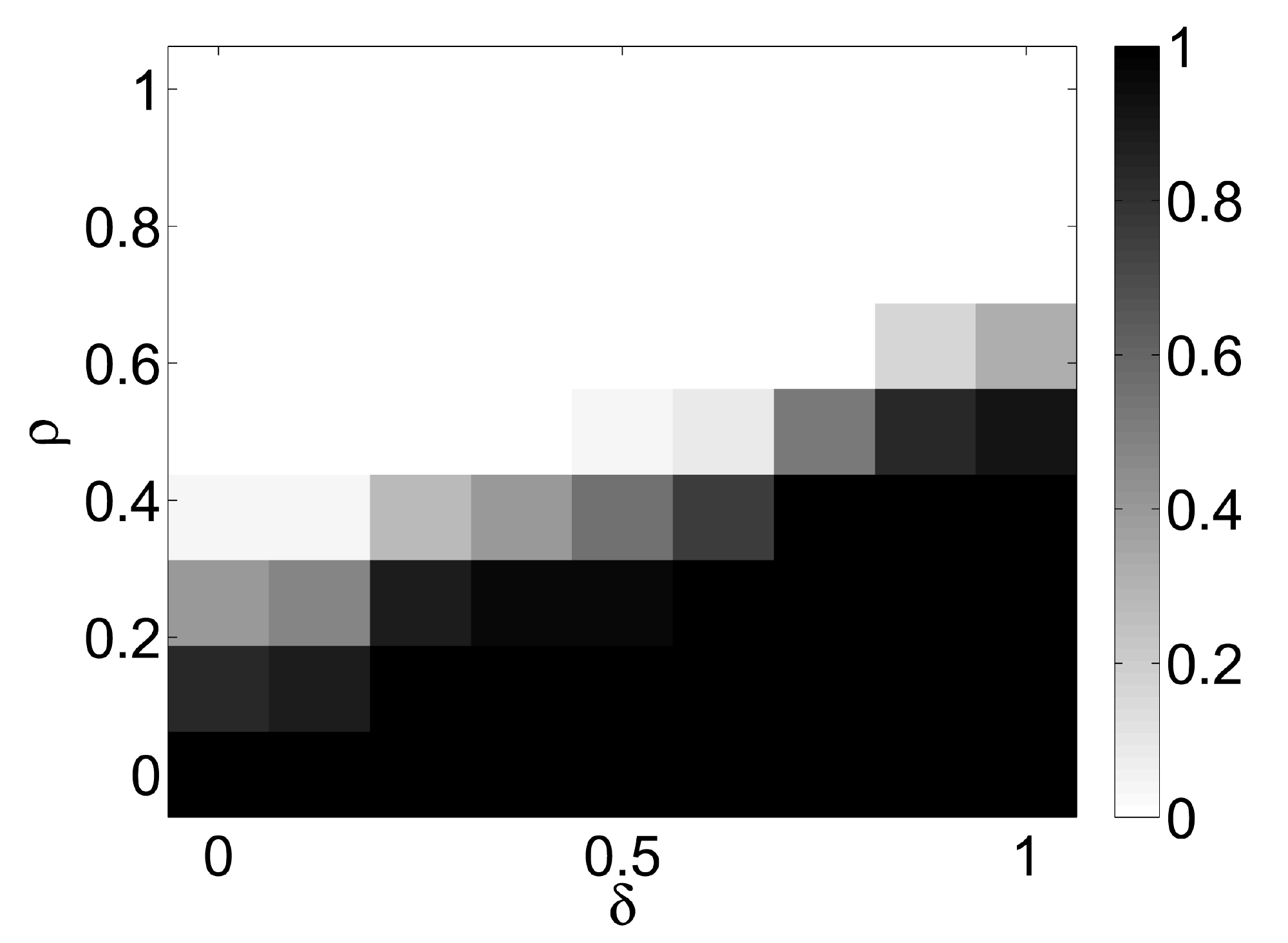} & (d) \includegraphics[width=1.6in]{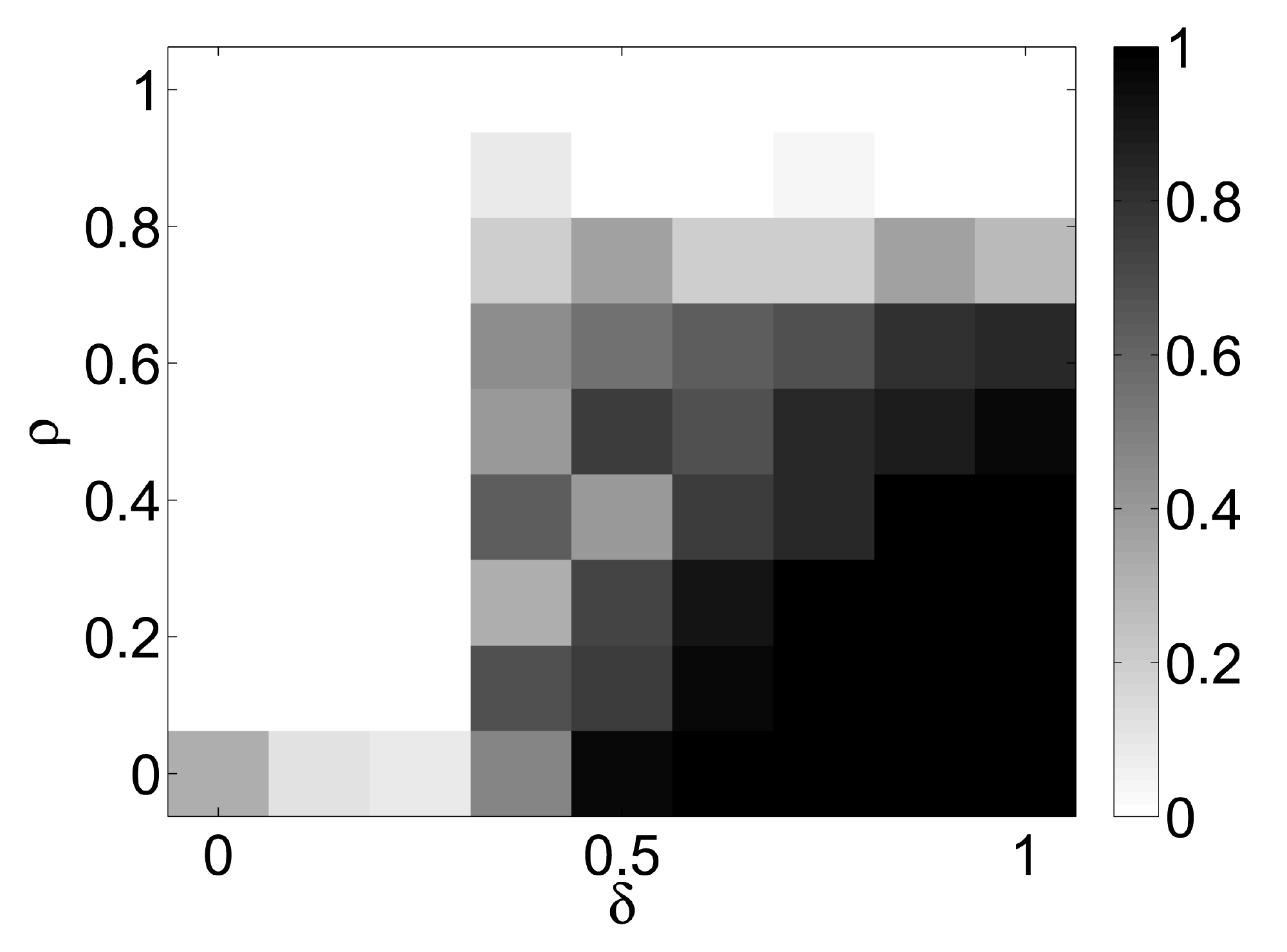}\\
    (b) \includegraphics[width=1.6in]{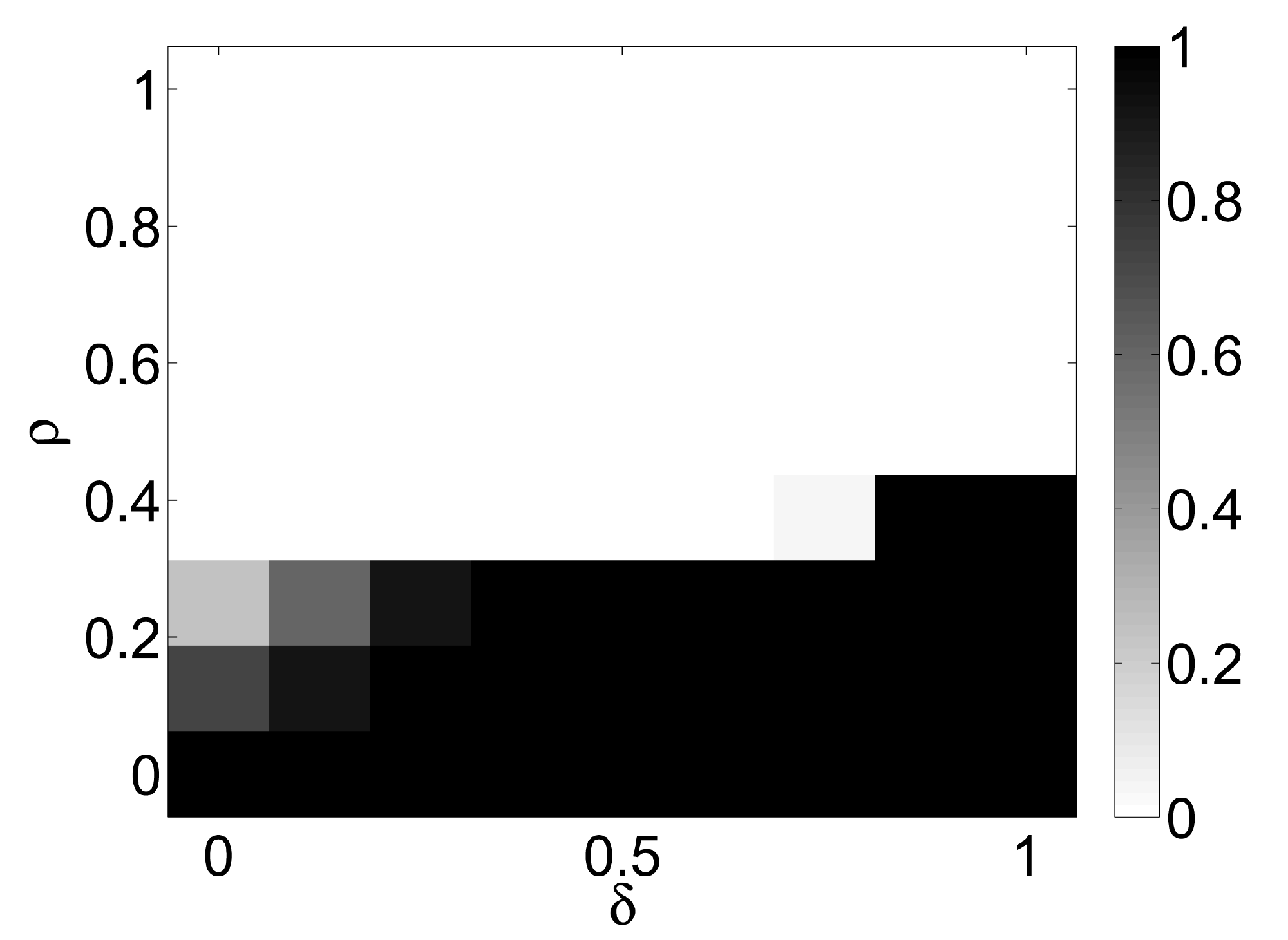} & (e) \includegraphics[width=1.6in]{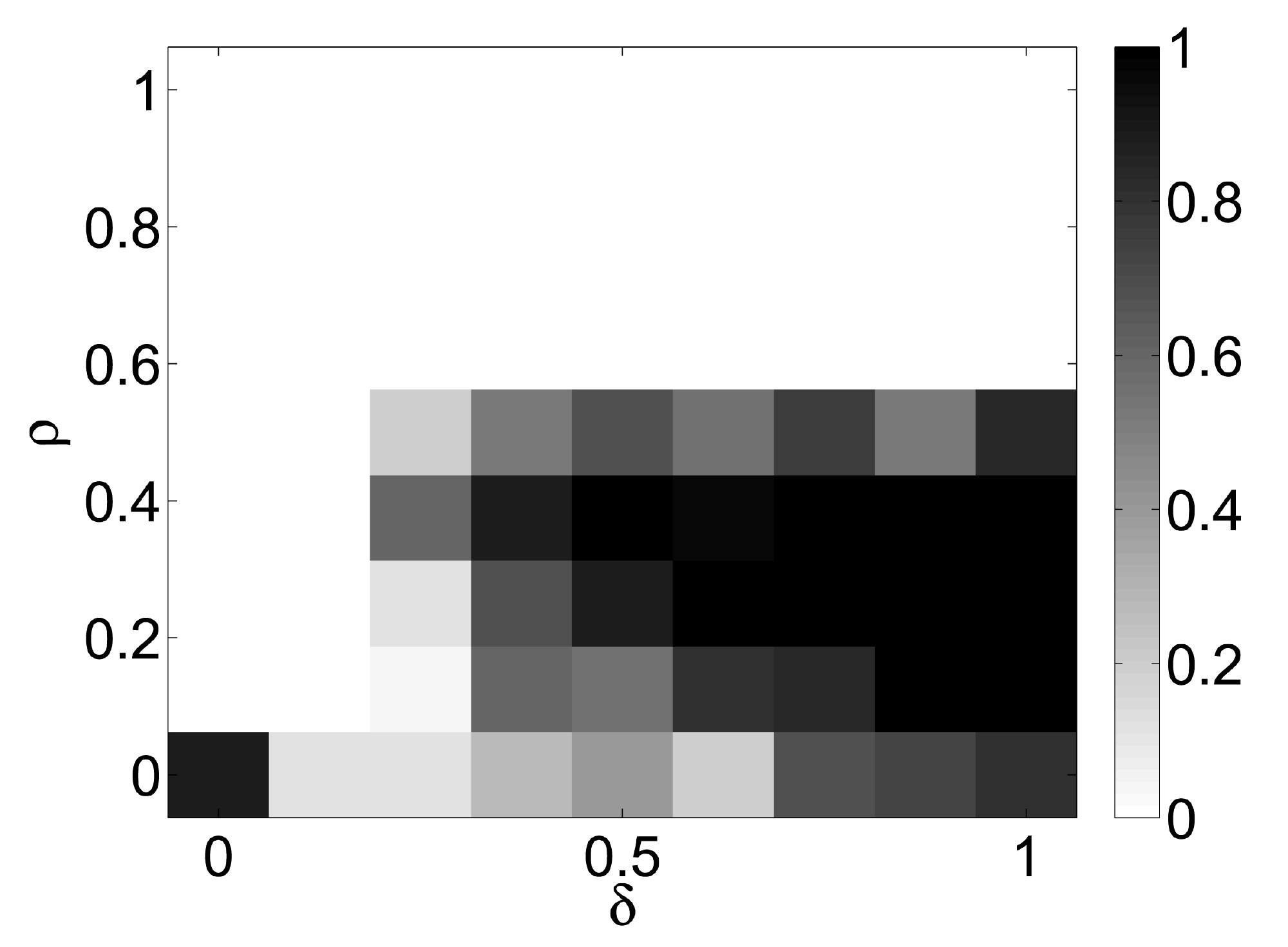} \\
    (c) \includegraphics[width=1.6in]{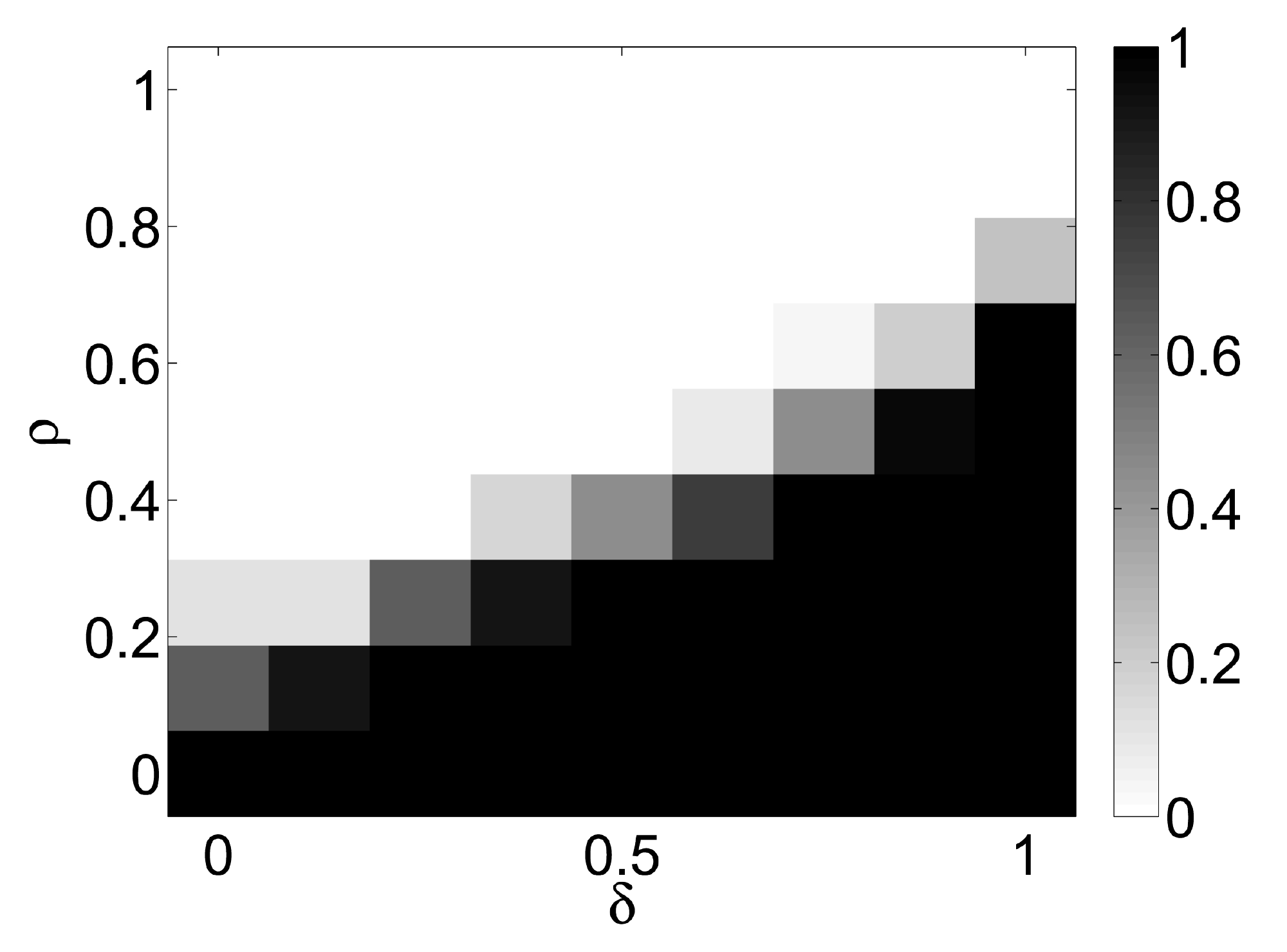} & (f) \includegraphics[width=1.6in]{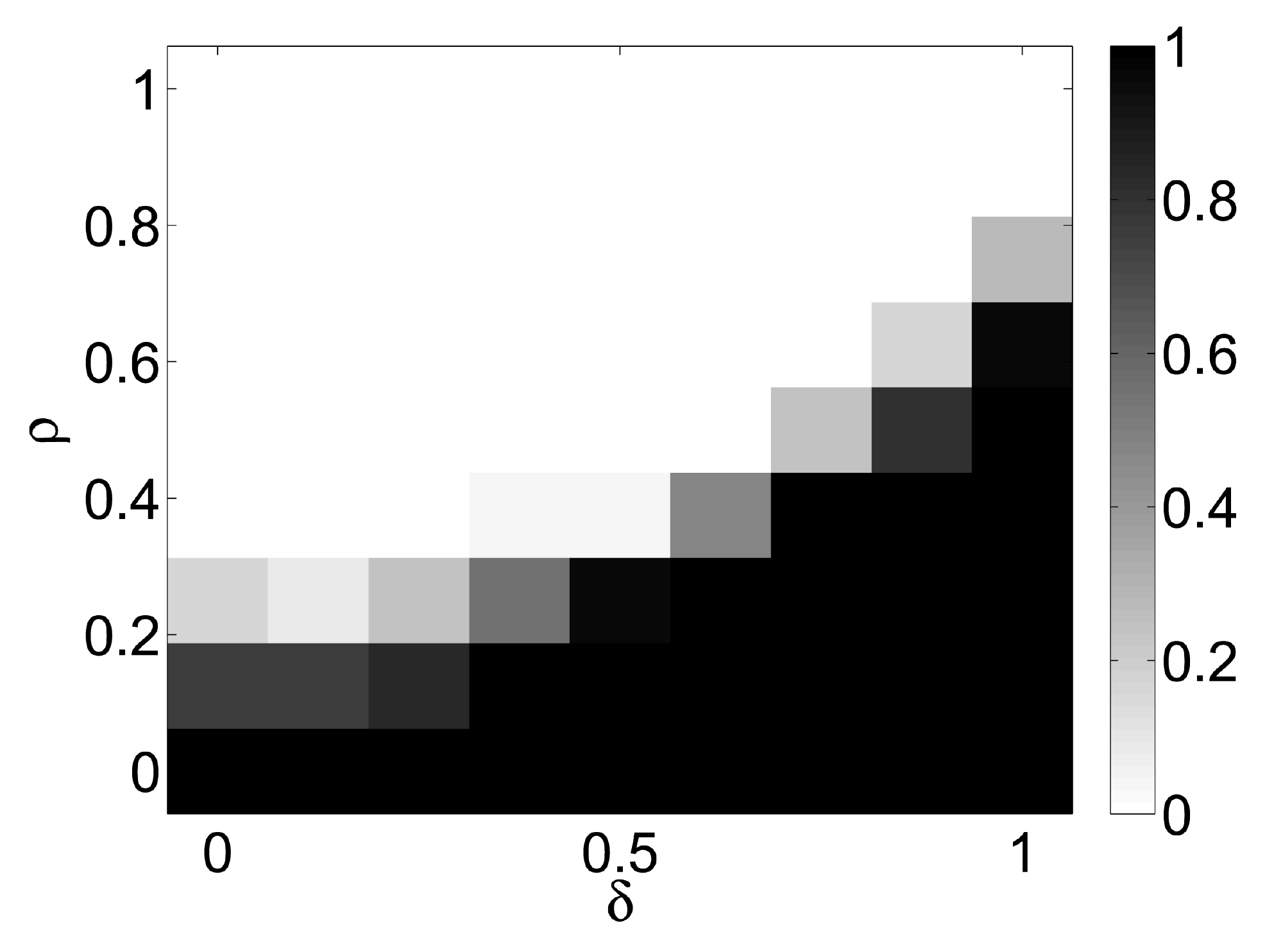} 
   \end{tabular}
   \vspace{-1mm}
   \caption{\small \sl Proportion of successes on Gaussian matrices using (a) PartInv, (b) CoSaMP and (c) $\ell_1$-minimization, and proportion of successes on correlated column subset matrices using (d) PartInv, (e) CoSaMP and (f) $\ell_1$-minimization for various values of $\delta=\frac{M}{N} \in (0,1)$ (horizontal axis) and $\rho=\frac{K}{M} \in (0,1)$ (vertical axis).
   \label{fig:gauss}}
   \vspace{-1mm}
\end{figure}

We show the results in Fig. \ref{fig:gauss}.

In the second case, we construct $M$ by $N$ matrices with $N=256$ and variable $M$ and a block diagonal structure. The columns are divided into 16 column subsets. In each subset we set $M/16$ rows to $1$ and add Gaussian noise with zero mean and variance $0.0625$. In addition, to every element of the matrix we add noise drawn from a zero-mean normal distribution with variance $\frac{1}{M}$. This produces heavy intra-subset correlation and light correlation across subsets. We let the coefficient vector $c$ contain $K$ non-zero elements drawn from a $N(0,1)$ distribution. We select 4 of the 16 subsets at random and in each subset select $\frac{K}{4}$ of the indices to have nonzero values, again uniformly at random. If some of the nonzeros were left over, they are accommodated in a fifth subset. For PartInv we set $L=\max\{K,0.8M\}$.

\subsection{Sensitivity of Recovery Performance to the Size of the Selected Subset}

\begin{figure}[!ht]
\vspace{1mm}
   \centering
   \begin{tabular}{@{\hspace{-2mm}}c@{\hspace{-3mm}}c}
    (a) \includegraphics[width=2in]{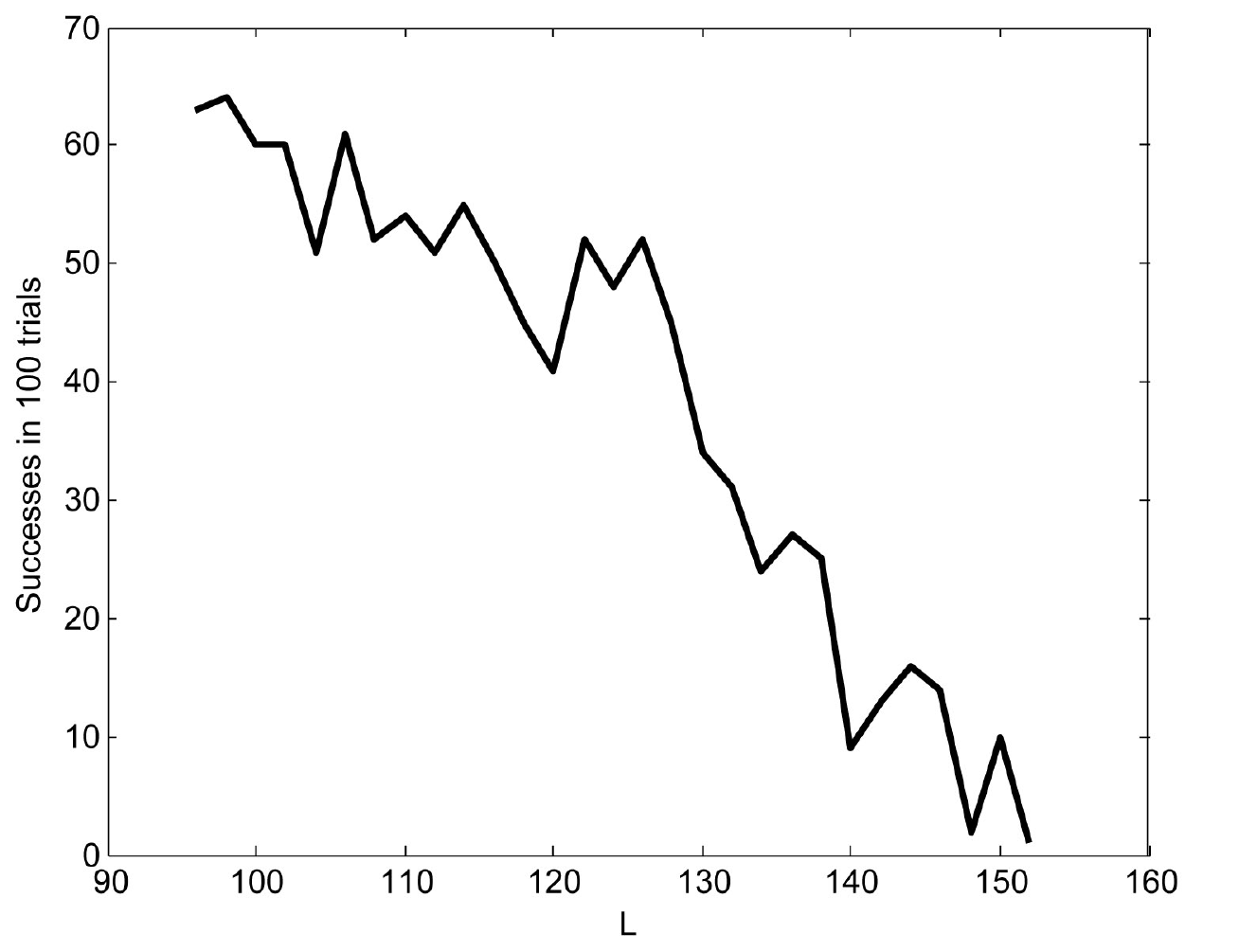} & (d) \includegraphics[width=2in]{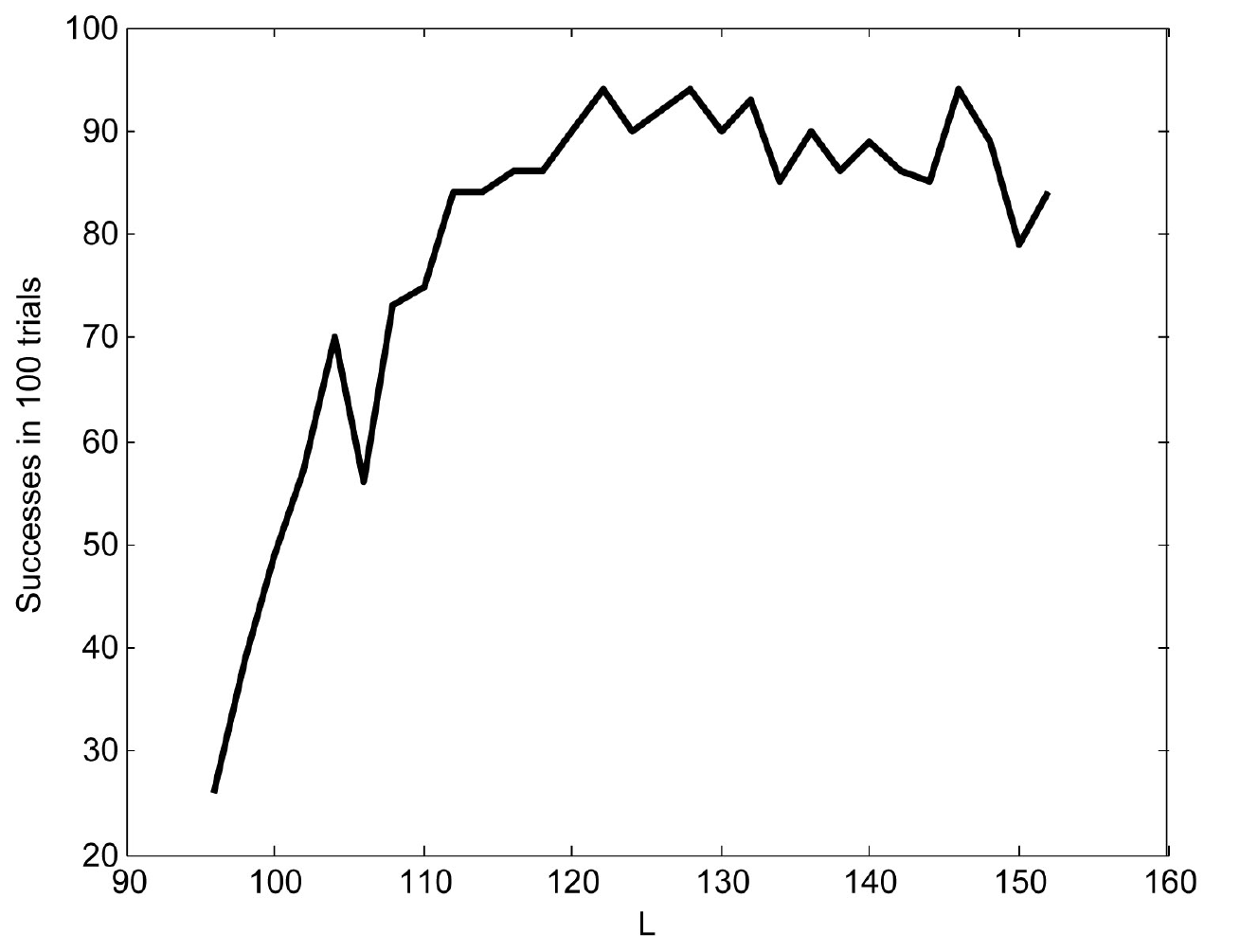}\\
    (b) \includegraphics[width=2in]{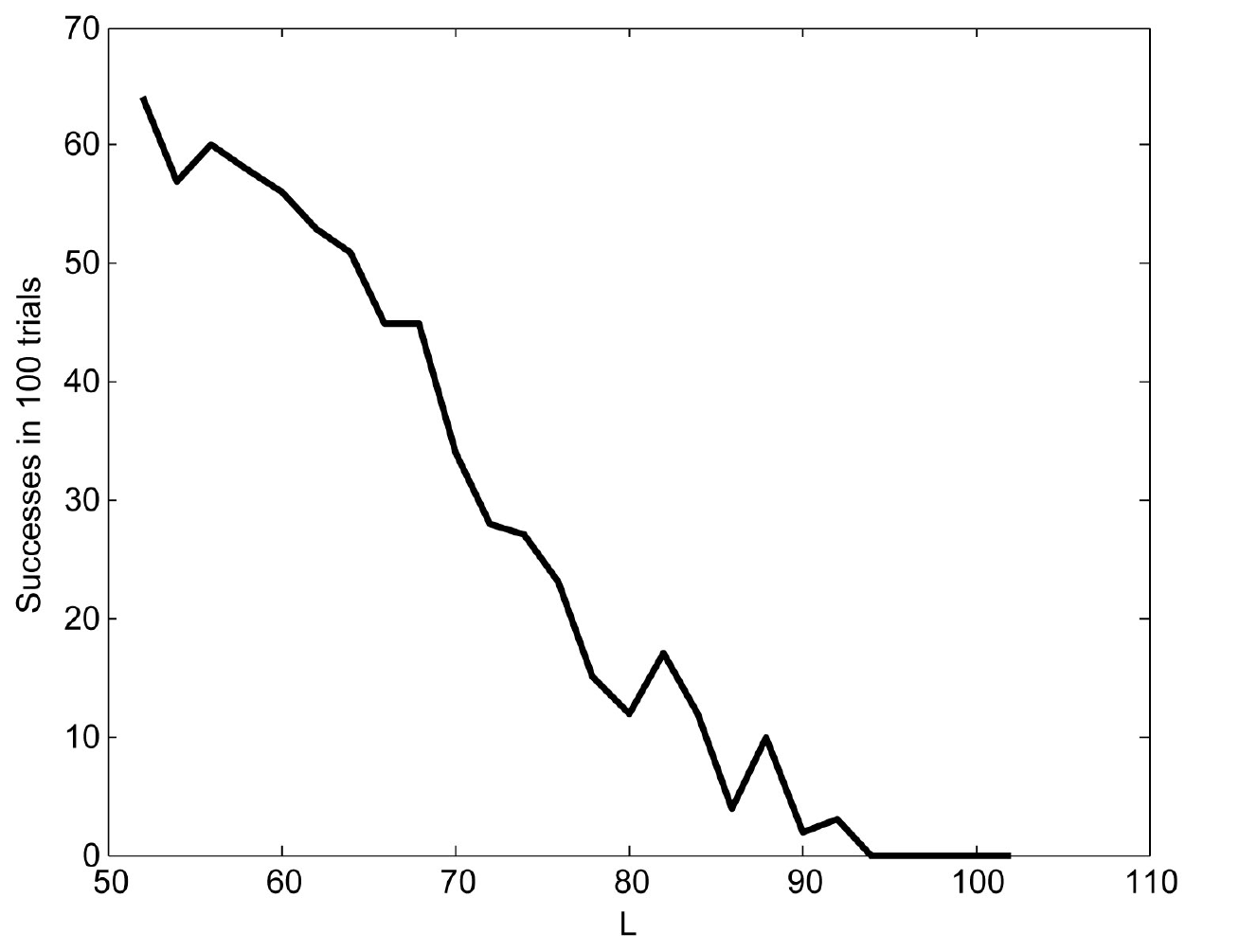} & (e) \includegraphics[width=2in]{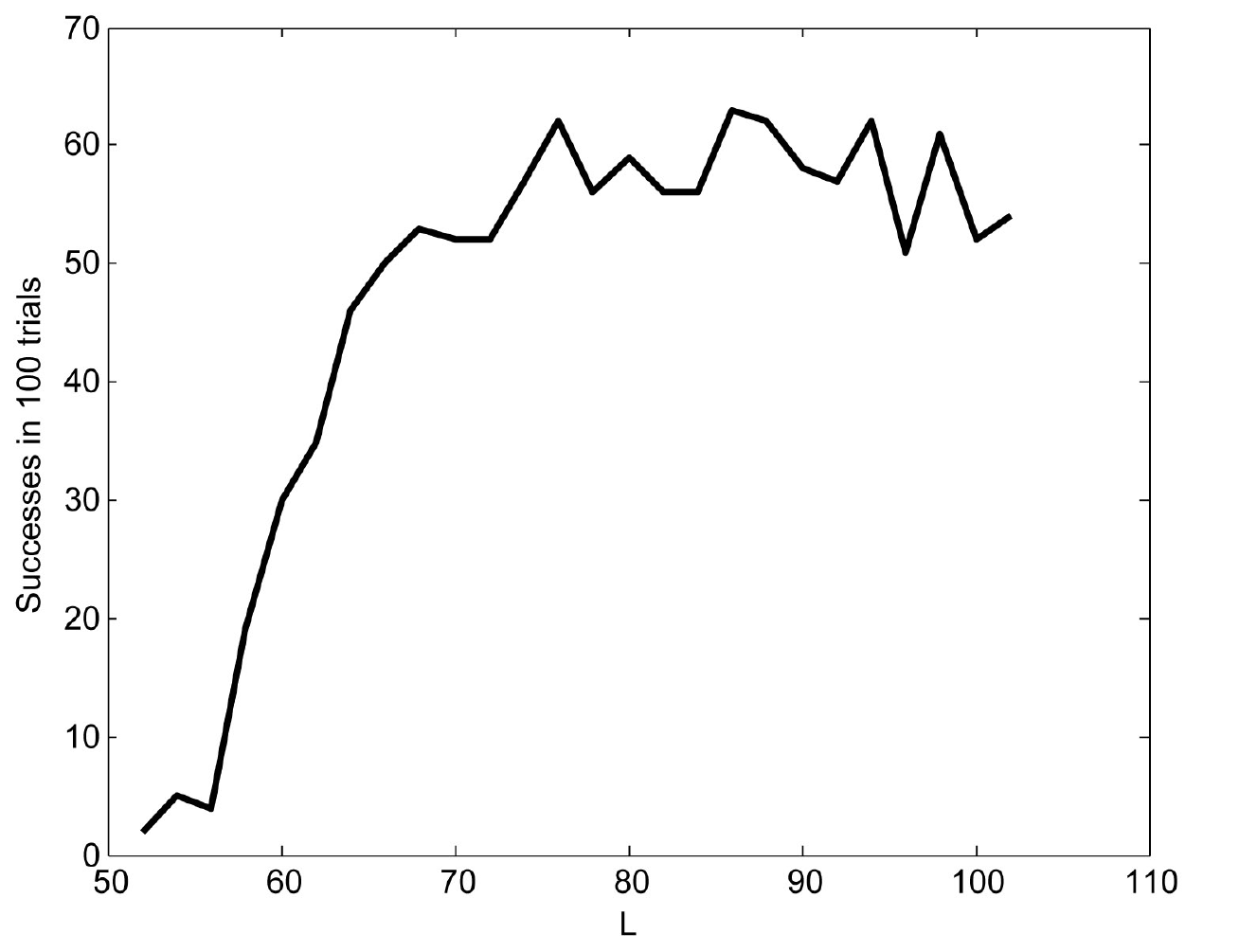} \\
    (c) \includegraphics[width=2in]{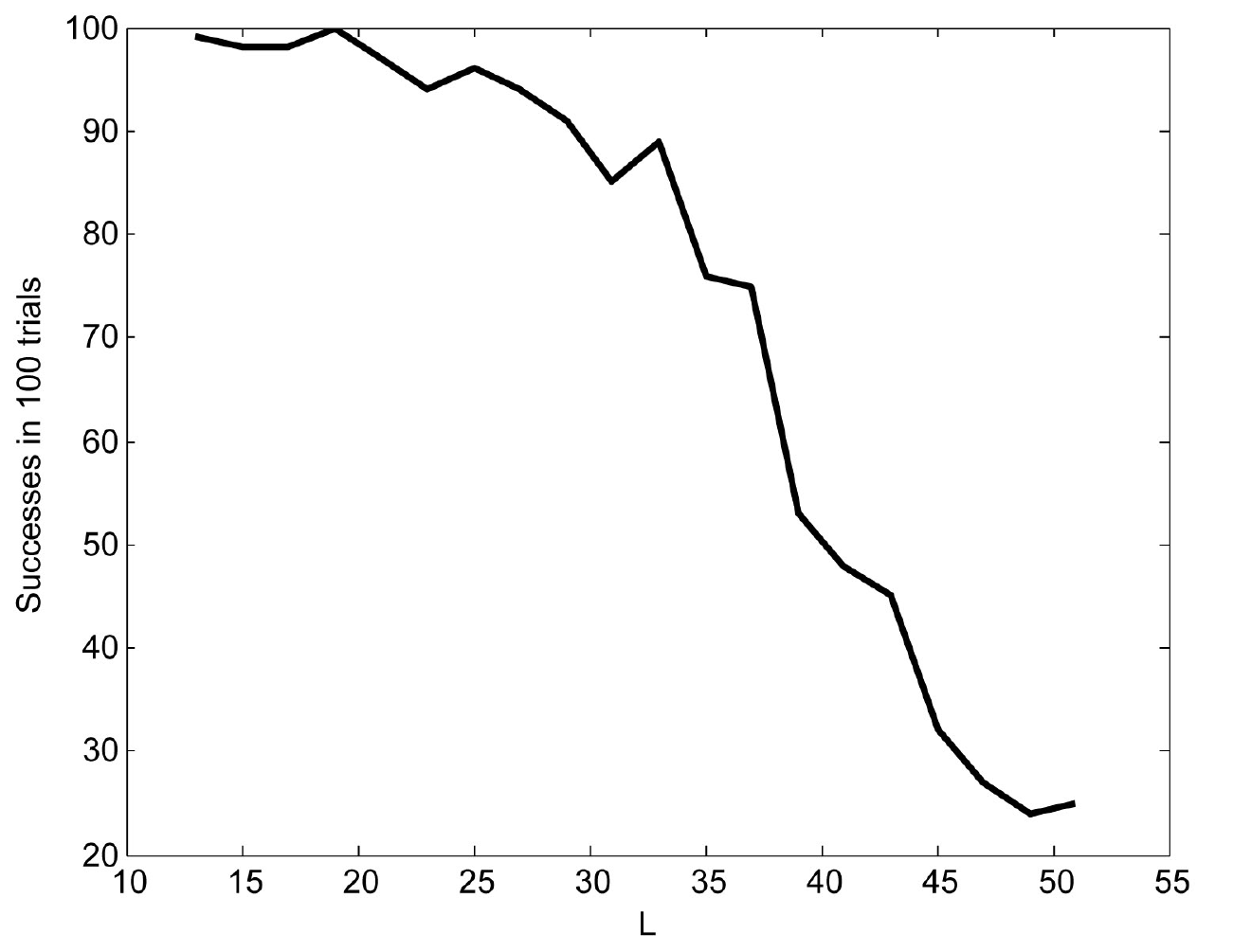} & (f) \includegraphics[width=2in]{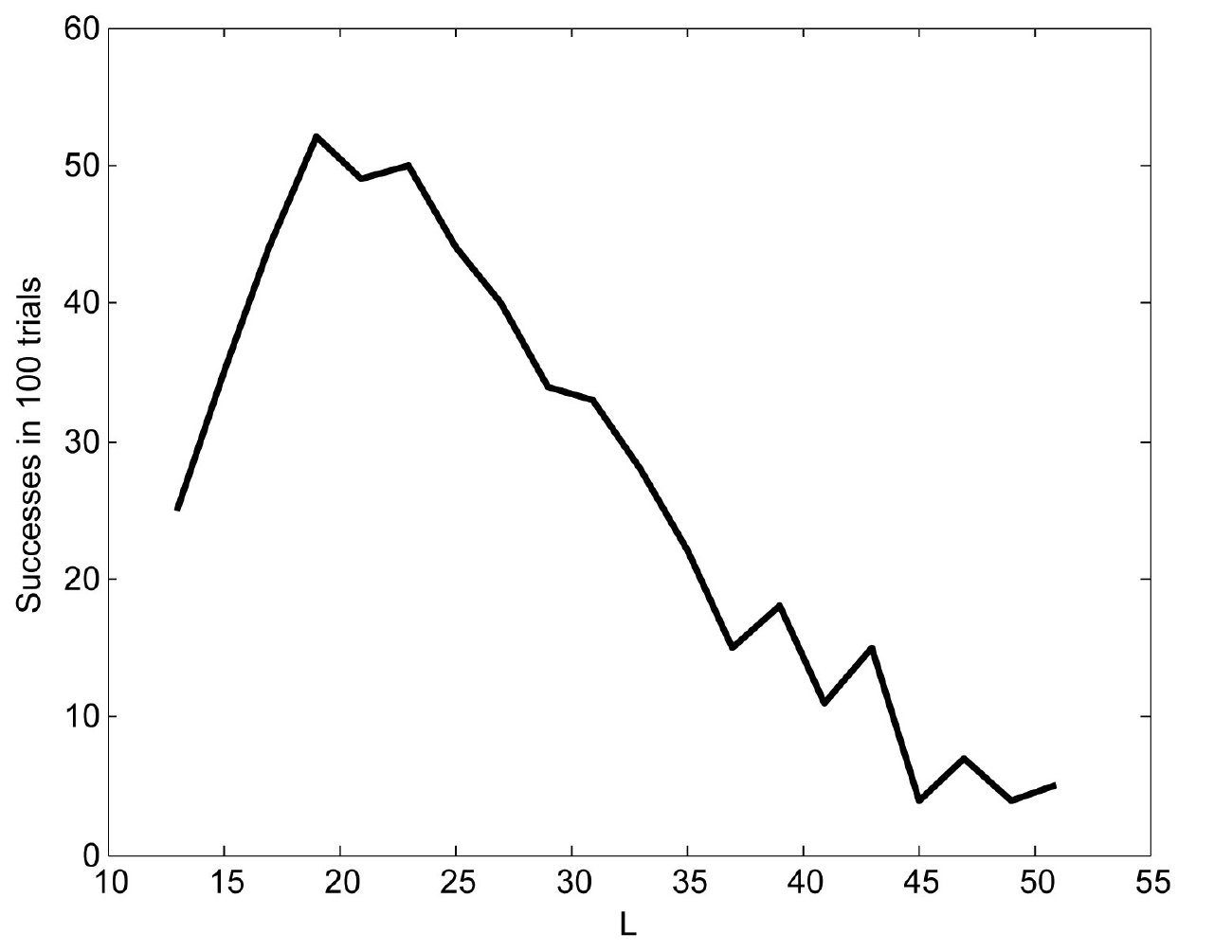} \\
   \end{tabular}
   \caption{\small \sl Proportion of successes using PartInv on Gaussian matrices(left) and correlated column subset matrices(right) for $M=192, K=96$ (a,d), for $M=128, K=52$ (b,e) and for $M=64, K=13$ (c,f) as L is varied from $K$ upto $0.8M$ in steps of 2.
   \label{fig:senseL}}
   \vspace{-1mm}
\end{figure}

We study the sensitivity of recovery performance to variations in parameter $L$, the size of the selected subset. We select three representative values of $M$ and $K$ that lie near the transition boundary between the high recovery and low recovery regions in the $\delta-\rho$ diagram (see Fig.~\ref{fig:gauss}). For each $(M,K)$ pair, we vary $L$ over the range $\{K..0.8M\}$ in steps of 2, for the Gaussian and correlated column subsets studied earlier. The results are shown in Fig.~\ref{fig:senseL}.

We see that for Gaussian matrices, the performance drops as $L$ increases from
$K$ towards $0.8M$, due to the presence of increasingly smaller singular values of $\Phi_{I}$ as $L$ approaches $M$, which lead to amplification of the noise component $\Phi_{I}^{\dagger}\Phi_{\tilde{I}}c_{\tilde{I}}$. For correlated column subset matrices, the performance is mixed, with performance improving as $L$ increases for the the first two cases, while it decreases for the last. Understanding this behavior is a topic for future work.

\begin{figure}[!ht]
\vspace{-1mm}
   \centering
   \begin{tabular}{@{\hspace{-2mm}}c@{\hspace{-3mm}}c}
   (a) \includegraphics[width=1.6in]{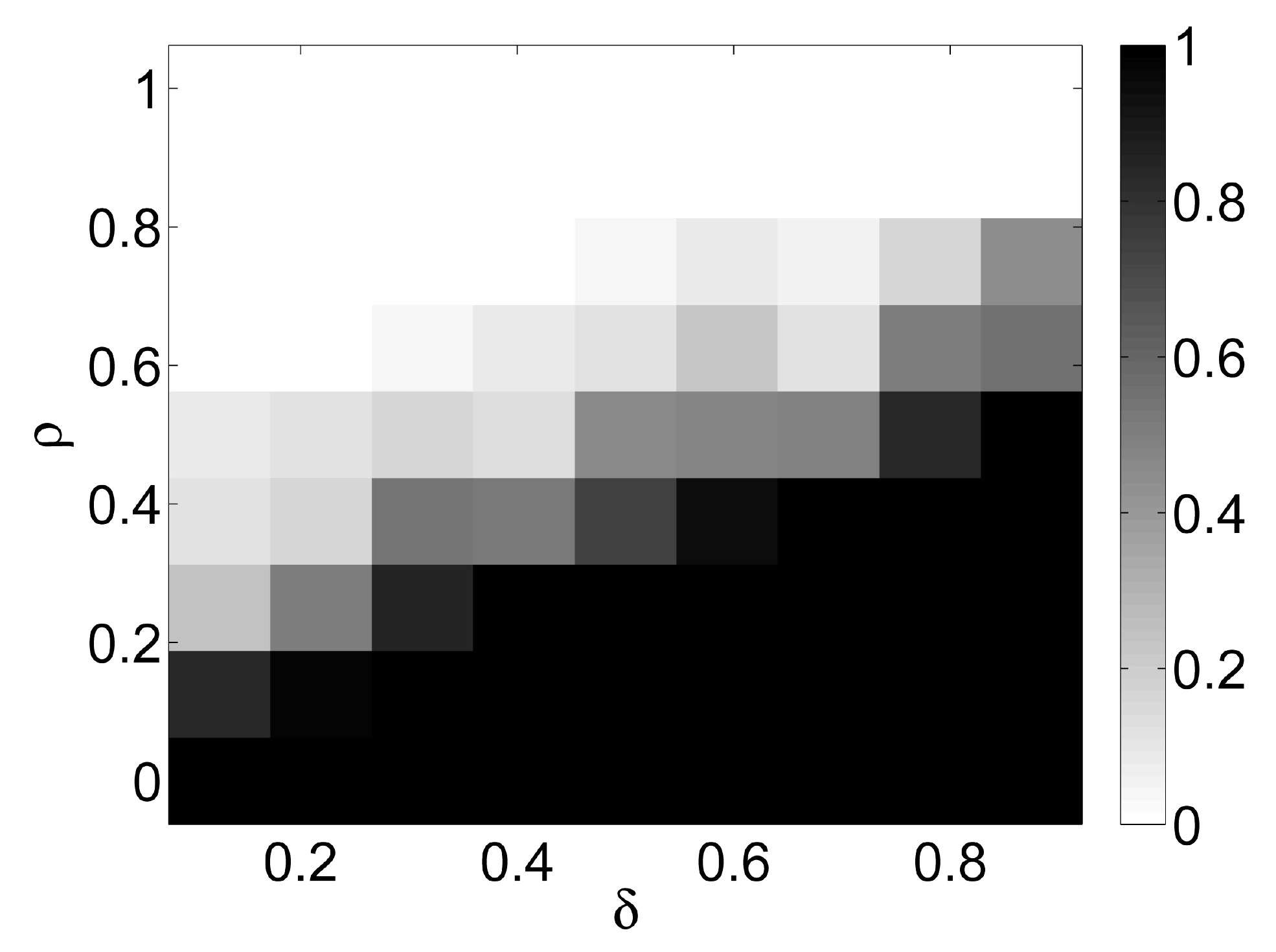} & (b) \includegraphics[width=1.6in]{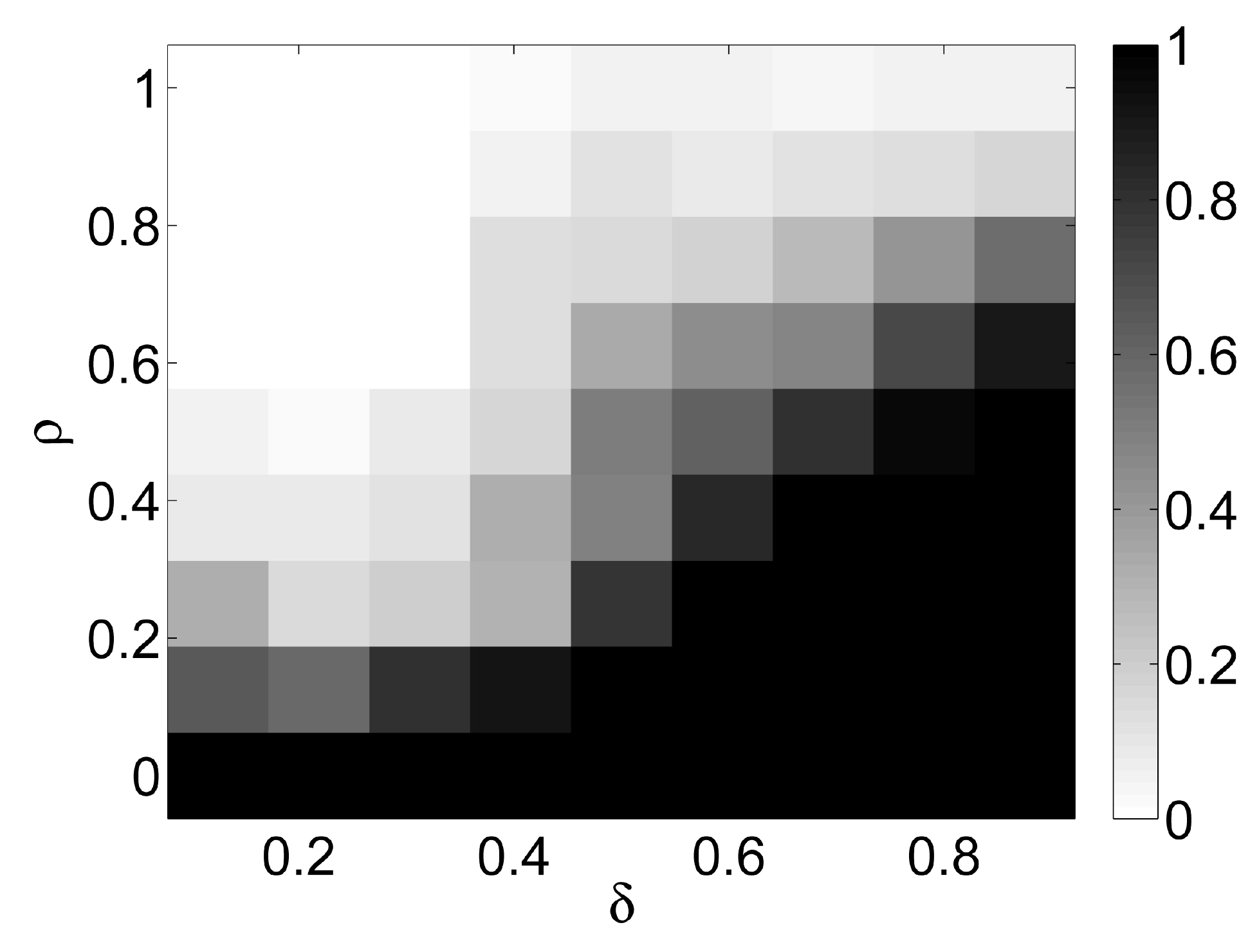}\\
   \end{tabular}
   \vspace{-1mm}
   \caption{\small \sl Best proportion of successes using PartInv as $L$ is varied from $K$ upto $M$ in steps of 2 for each $(M,K)$, on (a)Gaussian matrices (left) and (b) correlated column subset matrices(right).}
   \label{fig:bestPerf}
   \vspace{-1mm}
\end{figure}

We also obtain the best performance of Partial Inversion as $L$ is varied from $K$ to $M$ for each $\delta$ and $\rho$ in the range $(0,1)$ with 100 trials for each $(M,K,L)$ combination. The results are shown in Fig.~\ref{fig:bestPerf}. When compared with Fig.~\ref{fig:gauss} we see that there is a small improvement in the Gaussian case and a more significant improvement for the correlated columns case, especially for small $(\delta,\rho)$ values. This is consistent with the sensitivity to $L$ seen in Fig.~\ref{fig:senseL}, where we see that for correlated columns with $M=64,K=13$, the performance is best at lower values of $L$ close to $K$. We also provide the optimal values of $L$ for each $(\delta,\rho)$ in Table \ref{tab:bestL}.

\begin{table}[!ht]
\begin{center}
\begin{tabular}{|c|c|c|c|c|c|c|c|c|}
\hline
0  & 0  & 0  & 0  & 0  & 0  &  0  & 0 & 0 \\ \hline
0  & 0  & 0  & 0  & 0  & 0  &  0  & 0 & 0 \\ \hline
0  & 0  & 0  & 0  & 90 & 108 & 125 & 142 & 169 \\ \hline
0  & 0  & 54 & 61 & 76 & 91 & 109 & 128 & 162 \\ \hline
13 & 25 & 38 & 51 & 66 & 76  & 99 & 112 & 141 \\ \hline
10 & 20 & 30 & 44 & 51 & 65  & 77 & 83 & 92 \\ \hline
7  & 17 & 22 & 36 & 44 & 45  & 53 & 61 & 69 \\ \hline
5  & 10 & 17 & 20 & 25 & 30  & 35 & 40 & 46 \\ \hline
4  & 5  & 7  & 10 & 12 & 15  & 17 & 20 & 23 \\
\hline
\end{tabular}
\end{center}

\begin{center}
\begin{tabular}{|c|c|c|c|c|c|c|c|c|}
\hline
0 & 0 & 0 & 101 & 115 & 153 & 173 & 203 & 219 \\ \hline
0 & 0 & 0 &  99 & 126 & 130 & 177 & 199 & 206 \\ \hline
0 & 0 & 0 &  99 & 127 & 133 & 155 & 168 & 193 \\ \hline
0 & 0 & 0 &  85 & 116 & 129 & 133 & 164 & 180 \\ \hline
18 & 28 & 42 & 85 & 100 & 116 & 137 & 140 & 163 \\ \hline
16 & 20 & 76 & 100 & 87 & 93 & 131 & 109 & 112 \\ \hline
9  & 19 & 30 & 58  & 74 & 85 & 73  & 73 & 83 \\ \hline
9  & 14 & 21 & 50  & 69 & 46 & 49  & 48 & 54 \\ \hline
10 & 11 & 13 & 22  & 20 & 21 & 23  & 24 & 29 \\ 
\hline
\end{tabular}
\end{center}
\caption{\small \sl Values of $L$ that produce optimal performance for Gaussian(top) and correlated column subset matrices(bottom); a 0 indicates that the number of successes is zero} 
\label{tab:bestL}
\end{table}

\subsection{Recovery of Coefficients Concentrated on Wavelet Trees}

We next use Partial Inversion to recover nonzero coefficients that are concentrated on wavelet 
trees, which is commonly seen when a signal or image with discontinuities is decomposed in a wavelet basis. When the coefficients are concentrated on an isolated set (a set of columns that have low correlation with columns outside the set), a setwise estimator is especially useful
to identify the sets on which the coefficients are nonzero. 

Consider the 2D wavelet case.
Suppose that $I$ is the index set of columns of the wavelet basis belonging to a particular tree rooted at a coarse scale and containing finer scale coefficients. We have 
\begin{equation} 
\label{eq:z}
z_{I} = \Phi_{I}^{*}y 
      = \Phi_{I}^{*}\Phi_{I}c_{I} + \Phi_{I}^{*}\Phi_{\tilde{I}}c_{\tilde{I}}.
\end{equation}
Because $\Phi_{I}$ is relatively isolated from the columns in $\Phi_{\tilde{I}}$, the second term is small, and because most of the elements of $c_{I}$ are nonzero, the first term is large.
This is further intensified by the mutual correlation of the columns of $\Phi_{I}$ which is high because of the spatial overlap of the support of the wavelet bases in the tree. This motivates a simple selection criterion for measuring the strength of the nonzero coefficients in each wavelet tree $I$: $s_{I}=\sum\limits_{j\in I}|z_{j}|$. We use this criterion along with PartInv to select wavelet trees that are known to be nonzero. We denote the number of subsets by SetNum.
This modified method is described by Algorithm~\ref{PartInv-Wavelet}.

\begin{algorithm}[htbp]
\caption{Given $y=\Phi c$, with $K$ nonzero coefficients concentrated on wavelet trees, return best $K$-sparse approximation $\hat{c}$ }
\label{PartInv-Wavelet}
\begin{algorithmic}[1]
\STATE $\hat{c} \leftarrow \Phi^{*}y$; 
\STATE $k \leftarrow -1$
\FOR{$j=1 \to \text{SetNum}$}
\STATE $s_{j} \leftarrow \sum\limits_{l\in I_{j}}|\hat{c}_{l}|$
\ENDFOR
\STATE $I^{k+1} \leftarrow$ indices of columns contained in the sets with the largest magnitude $s_{k}$, to include at least $K$ coefficients.
\STATE $ k \leftarrow k+1$
\WHILE{Stopping condition not met} 
\STATE $\hat{c}_{I^{(k)}} \leftarrow \Phi_{I^{(k)}}^{\dagger} y$
\STATE $r \leftarrow y - \Phi_{I^{(k)}}\hat{c}_{I^{(k)}}$
\STATE $J^{(k)} \leftarrow \widetilde{I^{(k)}}$ 
\STATE $\hat{c}_{J^{(k)}} \leftarrow \Phi_{J^{(k)}}^{*}r$
\STATE Repeat lines $3-6$
\STATE $ k \leftarrow k+1$
\ENDWHILE
\end{algorithmic}
\end{algorithm}

\subsection{Experimental Results using wavelet trees}

To test this algorithm, we use the Daubechies-5 wavelet basis in two dimensions over $32\times 32$ size patches with 5 levels of decomposition. This gives a size $1024$ by $1024$ matrix $\Psi$. We divide this matrix into $49$ sets: $1$ set of the coarsest scale coefficients in a block of size $4\times 4$ containing the two coarsest scales, and $48$ other sets rooted at the coefficients at the next finer scale. Each of these sets contains $1+4+16=21$ coefficients in a quadtree structure. To create matrix $\Phi$ we first apply a blurring filter $H$ with a symmetric $5\times 5$ kernel that is close to a delta function. This simulates practical optical sample acquisition effects such as diffraction and helps prevent rank deficiency problems when carrying out inversion. We use different 2D sampling patterns to carry out the subsampling operation represented by matrix S. Hence the acquisition process is represented by $y=\Phi c$ where $\Phi=SH\Psi$. 
The sampling patterns are shown in Table~\ref{Table:Patterns} for each sampling rate $\delta=\frac{M}{N}$ used to generate the results. Each pattern is replicated 8 times in horizontal and vertical directions to give the $32\times 32$ sampling pattern used for matrix $S$.  The sampling patterns were selected to allow $\delta$ to increase in constant steps over the whole range $(0,1)$, while distributing the samples as uniformly and symmetrically as possible.
Common subsampling patterns used in superresolution problems would have sampled only a part of this interval, and are likely to give similar results as the patterns used here that are closest to them in density.  The filter kernel is a $5\times 5$ kernel with $0.29$ at the center and $0.02$ in other locations. The signals are generated by uniformly selecting at random wavelet trees to make the sparsity of the signal the specified value.  The coefficients in these trees are set to values chosen from a standard normal distribution, and the rest are set to zero.

\begin{table}[htbp]
\centering
\subfloat[$\delta=\frac{2}{16}$] {
\begin{tabular}{|c|c|c|c|}
\hline
0 & 0 & 0 & 0 \\ \hline
0 & 1 & 0 & 0 \\ \hline
0 & 0 & 0 & 0 \\ \hline
0 & 0 & 0 & 1 \\ \hline
\end{tabular}}
\qquad
\subfloat[$\delta=\frac{4}{16}$] {
\begin{tabular}{|c|c|c|c|}
\hline
1 & 0 & 0 & 0 \\ \hline
0 & 0 & 1 & 0 \\ \hline
0 & 1 & 0 & 0 \\ \hline
0 & 0 & 0 & 1 \\ \hline
\end{tabular}}
\qquad
\subfloat[$\delta=\frac{6}{16}$] {
\begin{tabular}{|c|c|c|c|}
\hline
1 & 0 & 1 & 0 \\ \hline
0 & 1 & 0 & 1 \\ \hline
1 & 0 & 0 & 0 \\ \hline
0 & 0 & 1 & 0 \\ \hline
\end{tabular}}
\qquad
\subfloat[$\delta=\frac{8}{16}$] {
\begin{tabular}{|c|c|c|c|}
\hline
1 & 0 & 1 & 0 \\ \hline
0 & 1 & 0 & 1 \\ \hline
1 & 0 & 1 & 0 \\ \hline
0 & 1 & 0 & 1 \\ \hline
\end{tabular}}
\qquad
\subfloat[$\delta=\frac{10}{16}$] {
\begin{tabular}{|c|c|c|c|}
\hline
0 & 1 & 0 & 1 \\ \hline
1 & 0 & 1 & 0 \\ \hline
0 & 1 & 1 & 1 \\ \hline
1 & 1 & 0 & 1 \\ \hline
\end{tabular}}
\qquad
\subfloat[$\delta=\frac{12}{16}$] {
\begin{tabular}{|c|c|c|c|}
\hline
1 & 1 & 0 & 1 \\ \hline
0 & 1 & 1 & 1 \\ \hline
1 & 1 & 1 & 0 \\ \hline
1 & 0 & 1 & 1 \\ \hline
\end{tabular}}
\qquad
\subfloat[$\delta=\frac{14}{16}$] {
\begin{tabular}{|c|c|c|c|}
\hline
1 & 1 & 1 & 1 \\ \hline
1 & 1 & 0 & 1 \\ \hline
1 & 1 & 1 & 1 \\ \hline
0 & 1 & 1 & 1 \\ \hline
\end{tabular}}
\caption{Sampling Patterns}
\label{Table:Patterns}
\vspace{-1mm}
\end{table}
%
%
\begin{figure}[ht]
\vspace{-1mm}
   \centering
   \begin{tabular}{@{\hspace{-2mm}}c@{\hspace{-3mm}}c}
  (a) \includegraphics[width=1.6in]{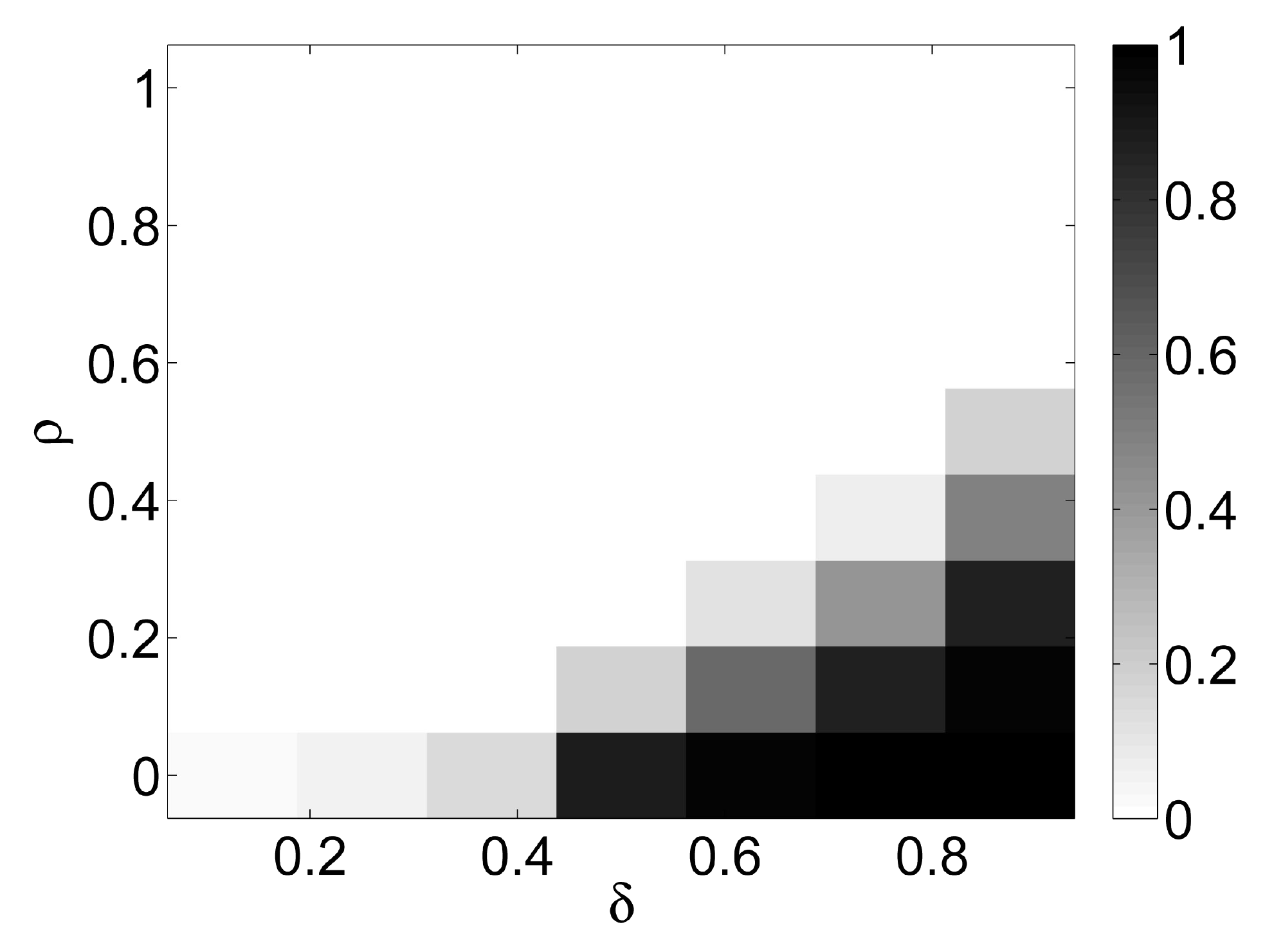} & (b) \includegraphics[width=1.6in]{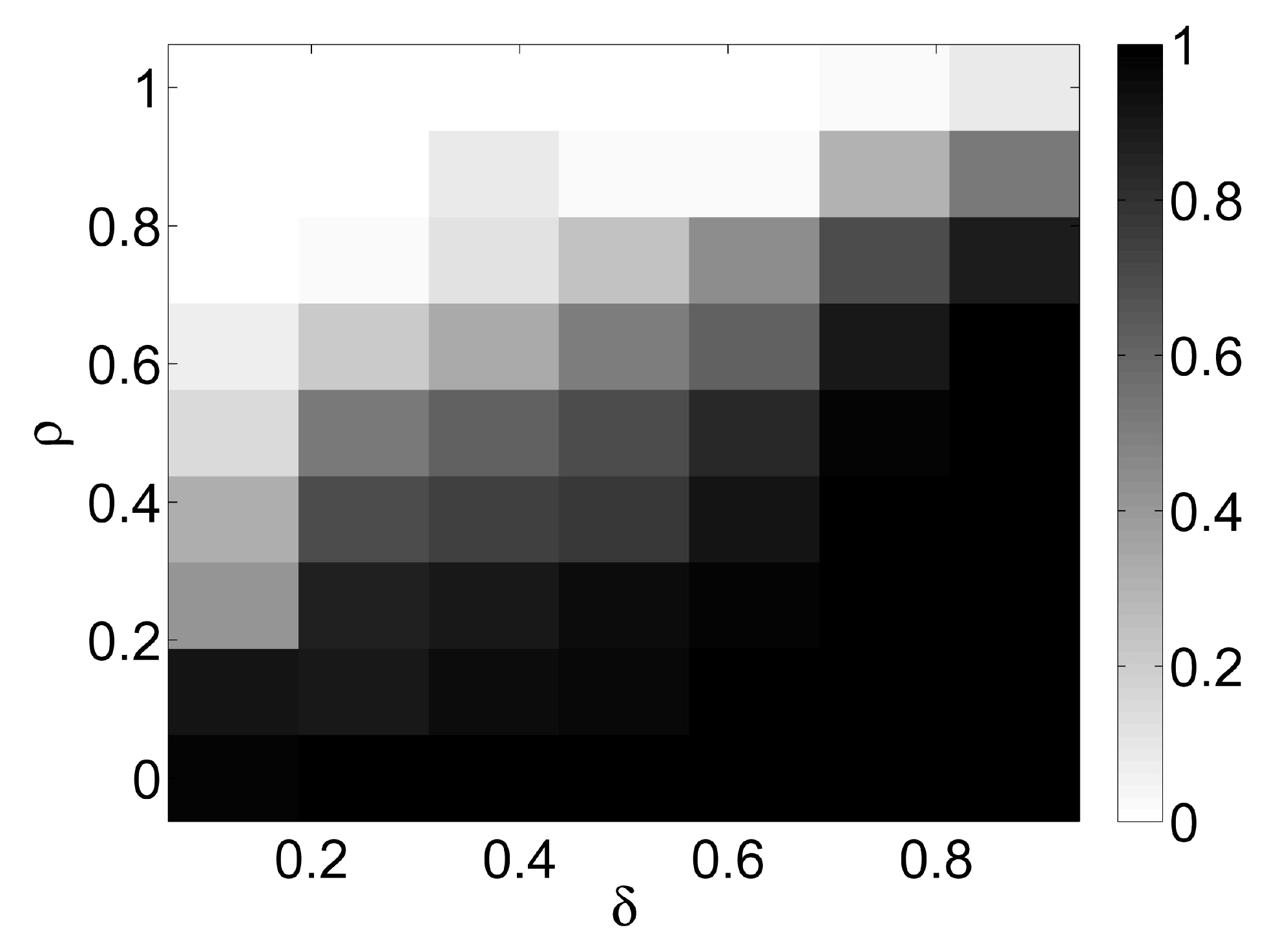}\\
   \end{tabular}
   \vspace{-1mm}
  \caption{\small \sl Proportion of successes with nonzero coefficients concentrated on wavelet trees from (a) $\ell_1$-minimization and (b) PartInv.
   \label{Table:Wavelet}}
   \vspace{-1mm}
\end{figure}

For each data point we carry out 100 trials.  We declare success if $\frac{1}{N}||c-\hat{c}||^{2}< 10^{-5}$ where $N=32\times 32$.  This shows improvement in selection performance with the sum estimator.  

\section{Conclusion and Future Work}\label{sec:conc}

We consider methods of compressive sensing recovery for sampling matrices  
with correlated columns.
This structure commonly arises in physical sample acquisition/reconstruction scenarios such as image super-resolution. We describe Partial Inversion, an algorithm that improves compressive sensing recovery by removing a source of noise in the initial estimator, and demonstrate its performance by simulations on Gaussian and correlated column subset matrices. We consider compressive sensing recovery when the nonzero coefficients are concentrated on wavelet trees, and demonstrate a simple estimator that improves selection of the trees that carry the nonzero coefficients.

We provide an analysis of the proposed method, which shows under mild conditions on the sensing matrix $\Phi$ that exact recovery is provided in the noiseless regime.  Our experimental results suggest that the method is also robust to noise, and that the assumptions placed on the correlation in the matrix can be weakened.  We believe an extension of our analysis is a good direction for future work. 

We also plan to consider compressive sensing recovery where the columns of the sampling matrix $\Phi$ can be grouped into nearly-isolated sets, such that correlation among pairs of columns within a set may be significant, but correlation between two columns that belong to different sets is relatively small. Future research will exploit this structure to efficiently reconstruct the signal.

\appendix

\section{Proof of Theorem~\ref{main}}
We start by pointing out that the first three $\Phi$ assumptions in Theorem~\ref{main} imply the following bounds:
\begin{align*}
\|\Phi_I^* \Phi_{I^c \cap T } \| & \le \delta, \quad   \forall \ |I| \leq L \\
\|\Phi_I^\dag \Phi_{I^c \cap T } \| & \le \delta, \quad   \forall \ |I|\leq L.
\end{align*}
To see this, write $\Phi_I = U\Sigma V^*$, the SVD decomposition in reduced form (i.e., the diagonals of $\Sigma$ are all strictly positive). Then the third condition becomes  
\[ \|\Phi_I \Phi_I^\dag \Phi_{I^c \cap T }\| = \|UU^* \Phi_{I^c \cap T }\| =\|U^* \Phi_{I^c \cap T }\| \le {\delta}/{A},\]
which holds because $U,V$ are both orthogonal matrices. Therefore,
\begin{align*}
\|\Phi_I^* \Phi_{I^c \cap T }\| & = \|V\Sigma U^*\Phi_{I^c \cap T}\| = \|\Sigma U^*\Phi_{I^c \cap T}\| \\
				     &\le \|\Sigma\| \cdot \|U^*\Phi_{I^c \cap T}\| = \|\Phi_I\| \cdot \|U^*\Phi_{I^c \cap T}\| \\
				     &\le A \cdot \delta/A = \delta. \\
\|\Phi_I^\dag \Phi_{I^c \cap T }\| & = \|V\Sigma^{-1} U^*\Phi_{I^c \cap T}\| = \|\Sigma^{-1} U^*\Phi_{I^c \cap T}\| \\
				     &\le \|\Sigma^{-1}\| \cdot \|U^*\Phi_{I^c \cap T}\| = \|\Phi_I^\dag\| \cdot \|U^*\Phi_{I^c \cap T}\| \\
				     &\le A \cdot \delta/A = \delta.
\end{align*}

We fix an iteration $k$ and consider $\chatIk$,  the approximation of $c$ at the beginning of iteration $k$. Let $T= \supp(c)$ be the true support of the signal $c$, so that $c=c_T$. Denote $\Jk = \widetilde{\Ik}$, the complement of $\Ik$ in the set of all column indices $1:N$. We then have
\begin{align}
\chatIk &= \PhiIkpinv y = \PhiIkpinv \Phi c = \PhiIkpinv (\PhiIk \cIk + \PhiJk \cJk) \nonumber \\
            &= \cIk + \PhiIkpinv \PhiJk \cJk. \label{eq:chatIk}
\end{align}
This equation is the foundation for the analysis below.
We consider the following three cases, depending on whether $T$ has been partially recovered, fully recovered, or not at all recovered.  Considering these cases separately serves to highlight how the algorithm gathers the support, and why each assumption is needed.  We aim to show that at each iteration, at least one new support element is identified, and that no correct support elements are lost.

\subsection*{Case 1: $T\subset \Ik$}
In this case we know that $\cJk = 0$ (which also implies that $y = \PhiIk \cIk$). So we obtain from \eqref{eq:chatIk} that $\chatIk = \cIk$, and hence, $r = y -\PhiIk \chatIk = 0$. This yields that $\chatJk = \PhistarJk r = 0$. As a result, the estimator of $c$ at iteration $k$ is $\chatk = c$. We have recovered the original signal $c$ without error.

In sum, if at one iteration $\Ik$ contains the entire support of the signal, then the corresponding estimator $\hat c$ coincides with $c$.

\subsection*{Case 2: $T\cap \Ik = \emptyset$}
In this case we have $\cIk=0$ and $T\subset \Jk$, so that $y=\PhiJk \cJk = \PhiT\cT$. 
Starting with \begin{equation*}\chatIk =\PhiIkpinv y,\end{equation*}
we get
\begin{align*}
r &= y-\PhiIk \chatIk = (I-\PhiIk \PhiIkpinv) y;\\
\chatJk &= \PhistarJk r =   \PhistarJk (I-\PhiIk \PhiIkpinv) y.
\end{align*}
Now plug in $y = \PhiT \cT$ to the last three equations:
\begin{align*}
\chatIk &= \PhiIkpinv  \PhiT \cT   \\
r           &= (I-\PhiIk \PhiIkpinv) \Phi_T c_T \\
\chatJk &= \PhistarJk (I-\PhiIk \PhiIkpinv) \PhiT \cT \nonumber \\
            &= \PhistarJk \PhiT \cT - \PhistarJk \PhiIk \PhiIkpinv \PhiT \cT
\end{align*}
Since $T\subset \Jk$, we may split $\chatJk$ into $\chatT$ and $\chatJkminusT$ and handle them separately:
\begin{align*}
\chatT&= \PhistarT \PhiT \cT - \PhistarT \PhiIk \PhiIkpinv \PhiT \cT; \\
\chatJkminusT&= \PhistarJkminusT \PhiT \cT - \PhistarJkminusT \PhiIk \PhiIkpinv \PhiT \cT.
\end{align*}
Using the assumptions we have the following estimates
\begin{align*}
\|\chatIk \|_2 
		&=  \|\PhiIkpinv  \PhiT \cT \|_2 \le \delta \|\cT\|_2 = \delta \|c\|_2 \\
\|\hat{c}_T\|_2   
	           &= \|\PhistarT\PhiT \cT - \PhistarT\PhiIk\PhiIkpinv \PhiT \cT \|_2 \\
		&\ge\|\PhistarT\PhiT \cT \|_2- \|\PhistarT\PhiIk\PhiIkpinv \PhiT \cT \|_2 \\
		&\ge (1-\delta)^2 \|c_T\|_2 - \delta  \| \PhiIkpinv \PhiT \cT \|_2 \\
		&\ge (1-\delta)^2 \|c_T\|_2 - \delta \cdot \delta \| \cT \|_2 \\
		&= (1-2\delta) \| c\|_2
\end{align*}
These two inequalities imply that the largest (in magnitude, same below) possible element in $\chatIk$ is $\delta \|c\|$, while there is at least one element in $\chatT$ that exceeds $\frac{1-2\delta}{\sqrt{K}} \|c\|$, where $K=|T|$ is the sparsity of the signal $c$.
Regarding the elements in $\chatJkminusT$, we upper bound them individually:  For any $i\in \JkminusT$, we have
\begin{align*}
|\hat c_{i}|  
		& = |\Phi^*_{i}\PhiT \cT - \Phi^*_{i}\PhiIk\PhiIkpinv \PhiT \cT | \\
		&\le|\Phi^*_{i}\PhiT \cT | + |\Phi^*_{i}\PhiIk\PhiIkpinv \PhiT \cT |  \\
		&\le \delta \|\cT\|_2 + \|\Phi_{i}\|_2 \cdot \|\PhiIk\PhiIkpinv \PhiT \cT\|_2 \\
		&\le \delta \|\cT\|_2 + 1\cdot \frac{\delta}{A}\|\cT\|_2 \\
		& \le 2\delta \|c\|_2.
\end{align*}
That is, the largest possible element in $\chatJkminusT$ is $2\delta \|c\|$.
Therefore, if we have that
\begin{equation}
\frac{1-2\delta} {\sqrt{K}}> 2\delta, \qquad \textrm{or} \quad \delta <  \frac{1}{2\sqrt K + 2}
\end{equation}
then at least one index from $T$ will be selected at the end of this iteration. Indeed, this inequality is true (provided that $K>3$) because $\delta<\frac{1}{3\sqrt{K}}$.

To summarize, if the initial $\Ik$ satisfies $\IkcapT = \emptyset$, at least one index from $T$ will be selected at the end of the iteration.

\subsection*{Case 3: $\IkcapT \ne \emptyset$, and $\JkcapT \ne \emptyset$}
This case is in between the first two cases. Our goal is to show that at the end of iteration $k$, all indices in $\IkcapT$ will be preserved and at least one new index from $\JkcapT$ will be selected.

To establish this, first note that 
\begin{equation} 
y = \PhiT \cT =  \PhiIkcapT \cIkcapT +  \PhiJkcapT \cJkcapT .
\end{equation} 
We continue directly from Equation \eqref{eq:chatIk} as follows:   
\begin{align*}
\chatIk &= \cIk + \PhiIkpinv \PhiJk \cJk  \\ 
	 &= \cIkcapT +  \PhiIkpinv \PhiJkcapT \cJkcapT
\end{align*}
Under the same assumptions, we split $\chatIk$ as $\chatIkcapT + \chatIkminusT$ and estimate them separately. Since \begin{equation*}
\|\PhiIkpinv \PhiJkcapT \cJkcapT\|_2 \le \delta \| \cJkcapT\|_2,
\end{equation*}
we know that every element of $\chatIkminusT$ is at most $ \delta \| \cJkcapT\| $ and every element of $\chatIkcapT$ is at least 
\begin{equation*}
3\delta \|c\|_2 -  \delta \| \cJkcapT\|_2 >2 \delta \| \cJkcapT\|_2.
\end{equation*}
We derive a formula for $\chatJk $ and will discuss its two parts $\chatJkcapT + \chatJkminusT$ also separately.
\begin{align*}
r           &=  y-\PhiIk \chatIk \\
 	 &=y - \PhiIk  \cIk - \PhiIk\PhiIkpinv \PhiJk \cJk \\
	 & = \PhiJkcapT \cJkcapT -  \PhiIk\PhiIkpinv \PhiJkcapT \cJkcapT\\
	 & = (I-\PhiIk \PhiIkpinv) \PhiJkcapT \cJkcapT \\
\chatJk &= \PhistarJk r  \\
	&= \PhistarJk (I-\PhiIk \PhiIkpinv) \PhiJkcapT  \cJkcapT \\
            &= \PhistarJk \PhiJkcapT \cJkcapT - \PhistarJk \PhiIk \PhiIkpinv \PhiJkcapT \cJkcapT.
\end{align*}
From above and using similar arguments, we have
\begin{align*}
\|\chatJkcapT\|_2  &= \| \PhistarJkcapT \PhiJkcapT \cJkcapT - \PhistarJkcapT \PhiIk \PhiIkpinv \PhiJkcapT \cJkcapT\|_2 \\
                          & \ge (1-\delta)^2 \|\cJkcapT\|_2 - \delta^2 \|\cJkcapT\|_2 \\
                          & =(1-2\delta) \|\cJkcapT\|_2
\end{align*}
and each element of $\chatJkminusT$ is no more than 
\begin{equation}
\delta\cdot \|\cJkcapT\|_2 + 1\cdot \frac{\delta}{A}\cdot \|\cJkcapT\|_2 \le 2\delta\|\cJkcapT\|_2.
\end{equation}

Considering $\chatIkcapT, \chatIkminusT, \chatJkcapT, \chatJkminusT$ all together, in order for all elements in $\chatIkcapT$ and at least one element from $\chatJkcapT$ to be among the $L$ largest entries of $\hat c$ and hence selected at the end of the iteration, we only need 
\begin{align*}
\frac{1-2\delta}{\sqrt{|\JkcapT|}} \|\cJkcapT\|_2&>2\delta\|\cJkcapT\|_2 
\end{align*}
or equivalently,
\begin{align*}
\delta & < \frac{1}{2\sqrt{|\JkcapT|}+2} .
\end{align*}
This inequality is indeed true because we already have $\delta < \frac{1}{2\sqrt{K}+2}$ from Case 2 and $|\JkcapT| < K$.  

\bibliographystyle{IEEEtran}
\bibliography{cs_struct}

\end{document}